\documentclass[a4paper,11pt,final]{article}
\usepackage{amsmath, amsthm, amssymb}
\usepackage{latexsym}
\usepackage[latin1]{inputenc}
\usepackage{epsf,exscale,times}
\usepackage[T1]{fontenc}
\usepackage[dvips]{graphicx}
\usepackage{graphics,epsfig}
\usepackage{epsfig,rotate,color}
\usepackage{showkeys}

\newcommand{\SR}{\mathcal{S}(\mathbb{R})}
\newcommand{\LdeuxR}{L^{2}(\mathbb{R})}

\newcommand{\LL}{\mathcal{I}}
\newcommand{\LLL}{\mathcal{L}}
\newcommand{\F}{\mathcal{F}}
\newcommand{\R}{\mathbb{R}}
\newcommand{\C}{\mathbb{C}}
\newtheorem{theoreme}{Theorem}
\newtheorem{lemme}{Lemma}
\newtheorem{definition}{Definition}
\newtheorem{Prop}{Proposition}
\newtheorem{Coro}{Corollary}
\newtheorem{remark}{Remark}

\title{A non-monotone conservation law for dune morphodynamics.}
\author{Natha\"el Alibaud, Pascal Azerad and Damien Is\`ebe}

\date{\today}

\begin{document}
\maketitle

\begin{quote} \footnotesize
\noindent \textsc{Abstract.} We investigate a non-local non linear
conservation law, first introduced by A.C. Fowler to describe
morphodynamics of dunes, see \cite{Fow01, Fow02}. A remarkable
feature is the violation of the maximum principle, which allows
for erosion phenomenon. We prove well-posedness for initial data
in $L^2$ and give explicit counterexample for the maximum
principle. We also provide numerical simulations corroborating our
theoretical results.
\end{quote}


\vspace{5mm}




\section{Introduction}

We investigate the following Cauchy problem:
\begin{equation}
\begin{cases}
\partial_t u(t,x) + \partial_x\left(\frac{u^2}{2}\right) (t,x) + \LL [u(t,.)] (x)
- \partial_{xx}^2 u(t,x) = 0 & t \in (0,T), x \in \mathbb{R}, \\
u(0,x)=u_0(x) & x \in \mathbb{R},
\end{cases}
\label{fowlereqn}
\end{equation}
where $T$ is any given positive time, $u_0 \in L^2(\R)$ and $\LL$
is a non-local operator defined as follows: for any Schwartz
function $\varphi \in \SR$ and any $x \in \mathbb{R}$,
\begin{equation*}
\LL [\varphi] (x) := \int_{0}^{+\infty} |\zeta|^{-\frac{1}{3}}
\varphi''(x-\zeta) d\zeta . \label{nonlocalterm}
\end{equation*}
\begin{remark}
Equation \eqref{fowlereqn} can also be written in conservative form
\begin{equation*}
\partial_t u + \partial_x \left(\frac{u^2}{2} + \LLL [u]-  \partial_x u\right) = 0
\end{equation*}
  where $$\LLL[\varphi] (x) :=
 \int_{0}^{+\infty} |\zeta|^{-\frac{1}{3}}
\varphi'(x-\zeta) d\zeta.$$
 \end{remark}
\bigskip

Equation \eqref{fowlereqn} appears in the work of Fowler
\cite{Fow01,Fow02} on the evolution of \textit{dunes}; the term
dunes refers to instabilities in landforms, which occur through
the interaction of a turbulent flow with an erodible substrate.
Equation \eqref{fowlereqn} is valid for a river flow (from left to
the right) over a erodible bottom $u(t,x)$ with slow variation.
For more details on the physical background, we refer the reader
to \cite{Fow01,Fow02}.

Roughly speaking, $\LL[u]$ is a weighted mean of second
derivatives of $u$ with the bad sign; hence, this term has a
deregularizing effect and the main consequence is probably the
fact that \eqref{fowlereqn} does not satisfy the maximum principle
(see below for more details). Nevertheless, one can see that the
diffusive operator $-\partial_{xx}^2$ controls the instabilities
produced by $\LL$ and ensures the existence and the uniqueness of
a smooth solution for positive times. The starting point to
establish these facts is the derivation of a new formula for the
operator $\LL$, namely \eqref{intformula}. This result allows
first to show easily that $\LL-\partial_{xx}^2$ is a
pseudo-differential operator with symbol $\psi_\LL(\xi)=4 \pi^2
\xi^2-a_\LL |\xi|^{\frac{4}{3}}+i \: b_\LL \xi
|\xi|^{\frac{1}{3}}$, where $a_\LL$ and $b_\LL$ are positive
constants (see \eqref{formule pseudo}). The symbol $4 \pi^2 \xi^2$
corresponds to the diffusive operator $-\partial^2_{xx}$ and
$-a_\LL |\xi|^{\frac{4}{3}}+i \: b_\LL \xi |\xi|^{\frac{1}{3}}$ is
the symbol of the nonlocal operator $\LL$. Notice that this last
symbol contains a fractional anti-diffusion
$-a_\LL|\xi|^{\frac{4}{3}}$ (recall that this is the symbol of
$-(-\partial_{xx}^2)^{\frac{4}{6}}$, up to a positive
multiplicative constant) and a fractional drift $i \: b_\LL \xi
|\xi|^{\frac{1}{3}}$. Because of the fact that the fractional
anti-diffusion is of order $\frac{4}{3}$, the real part of
$\psi_\LL(\xi)$ behaves as $\xi^2$, up to a positive
multiplicative constant, as $\xi \rightarrow +\infty$. A
consequence is that Equation \eqref{fowlereqn} has a regularizing
effect on the initial data: even if $u_0$ is only $L^2$, the
solution $u$ becomes $C^{\infty}$ for positive times. The
uniqueness of a $L^\infty((0,T);L^2)$ solution is obtained by the
use of a mild formulation (see Definition \ref{def Duhamel}) based
on  Duhamel's formula \eqref{formule duhamel}, in which appears
the kernel $K$ of $\LL-\partial_{xx}^2$. The use of such a formula
also allows to prove  local-in-time existence with the help of a
contracting fixed point theorem. Such an approach is quite
classical; we refer the reader, for instance, to the book of Pazy
\cite{Pazy} and the references therein on the application of the
theory of semigroups of linear operators to partial differential
equations. We also refer the reader to the work of Droniou
\textit{et al.} in \cite{DroGalVov03} for fractal conservation
laws of the form
\begin{equation}
\partial_t u+\partial_x (f(u))+(-\partial_{xx}^2)^{\frac{\lambda}{2}} [u]=0,
\label{drofrac}
\end{equation}
where $f$ is locally Lipschitz continuous and $\lambda \in (1,2]$,
and to the work of Tadmor \cite{Tad86} on the Kuramoto-Sivashinsky
equation:
$$
\partial_t u+\frac{1}{2}|\partial_x u|^2-\partial_{xx}^2 u=(-\partial_{xx}^2)^2 [u].
$$
In fact,  fractal conservation law \eqref{drofrac} is monotone and the
global existence of a $L^\infty$ solution is based on the fact
that the $L^\infty$ norm of $u$ does not increase. In our case,
this is not true and we have to use a classical energy estimate
 to get a global $L^2$ estimate. The
regularizing effect on the initial data are first proved by a
fixed point theorem on the Duhamel's formula to get $H^1$
regularity in space and next by a bootstrap method to get further
regularity. This technique has already been used in
\cite{DroGalVov03}.

On the other hand, one of our main result is probably the proof of
the failure of the maximum principle for \eqref{fowlereqn}: more
precisely, we exhibit positive dunes which take negative values in
finite time, since we establish that the bottom is eroded
downstream from the dune. We also give some numerical results that
illustrate this fact (for more precision, see Remark \ref{rem max}
and Section \ref{sect numerique}). The proof of the failure of the
maximum principle is based on the integral formula
\eqref{intformula}. Roughly speaking, this formula means that
$\LL$ is a L\'evy operator with a bad sign, see \cite{bertoin}.
Notice that the Kuramoto-Sivashinsky equation is also
non-monotone, but no proof of the failure of the maximum principle
is given in \cite{Tad86}.

The paper is organized as follows. In Section \ref{sect prelim},
we give the integral and pseudo-differential formula for
$\LL$; we also establish the properties on the kernel $K$ of
$\LL-\partial_{xx}^2$ that will be needed. In Section \ref{sect
duhamel}, we define the notion of mild solution for
\eqref{fowlereqn}. Sections \ref{sect uniqueness} and \ref{sect
existence} are, respectively, devoted to the proof of the
uniqueness and the existence of a mild solution; Section \ref{sect
existence} also contains the proof of the regularity of the
solution. The proof of the failure of the maximum principle is
given in Section \ref{sect princ max}. Finally, we give in Section
\ref{sect numerique} some numerical simulations that illustrate
the theory of the preceding sections.

\bigskip

Here are our main results.

\begin{theoreme}\label{theo principal}
Let $T>0$ and $u_0 \in \LdeuxR$. There exists a unique mild
solution $u \in L^\infty((0,T);L^2(\R))$ of \eqref{fowlereqn} (see
Definition \ref{def Duhamel}). Moreover,
\begin{itemize}
\item[i)] $u \in C^\infty((0,T] \times \R)$ and for all $t_0 \in
(0,T]$, $u$ and all its
derivatives belong to $C([t_0,T];L^2(\R))$.

\item[ii)] $u$ satisfies $\partial_t u + \partial_x(\frac{u^2}{2}) + \LL [u] -
\partial^2_{xx} u = 0$, on $(0,T] \times \R$, in the classical sense ($\LL
[u]$ being properly defined by \eqref{intformula} and \eqref{formule pseudo}).

\item[iii)] $u \in C([0,T];L^2(\R))$ and $u(0,.)=u_0$ almost everywhere (a.e. for short).
\end{itemize}
\end{theoreme}
\begin{Prop}[$L^2$-stability]
\label{L2stab}
Let
$(u,v)$ be solutions to \eqref{fowlereqn} with respective $L^2$
initial data $(u_0,v_0)$, we have:
$$
||u-v||_{C([0,T];L^2(\R))} \leq
C\left(T,M,||u_0||_{L^2(\R)},||v_0||_{L^2(\R)}\right)
||u_0-v_0||_{L^2(\R)}
$$
where $M : = \max \left( ||u||_{C([0,T];L^2(\R))},
||v||_{C([0,T];L^2(\R))}\right)$.
\end{Prop}

\begin{theoreme}[Failure of the maximum principle]\label{theo max}
Assume that $u_0 \in C^2(\R) \cap H^2(\R)$ is nonnegative and such
that there exist $x_\ast \in \R$ with $u_0(x_\ast)
=u_0'(x_\ast)=u_0''(x_\ast)=0$ and
\begin{equation*}
\int_{-\infty}^{0} \frac{u_0(x_\ast+z)}{|z|^{7/3}} \: dz>0.
\end{equation*}
Then, there exists $t_\ast>0$  with $u(t_\ast,x_\ast) < 0$.
\end{theoreme}

\begin{remark}\label{rem max}
Hypothesis of the theorem above are satisfied, for instance, for
non-negative $u_0 \in C^2(\R) \cap H^2(\R)$ such that there exists
$x_\ast \in \R$ with $u_0(x_\ast) =u_0'(x_\ast)=u_0''(x_\ast)=0$
and
$$
\forall x \leq x_\ast, \;  u_0(x)\geq 0
\quad \mbox{and} \quad
\exists x_0 <x_\ast \mbox{ s.t. } u_0(x_0)>0.
$$
A simple example of such an initial dune is shown in Figure
\ref{non_max_princ}. Observe that the bottom is eroded downstream
from the dune (recall that the nonlinear convective term
propagates a positive dune from the left to the right).
\begin{figure}[!ht!]
\begin{center}
\includegraphics[scale=0.2]{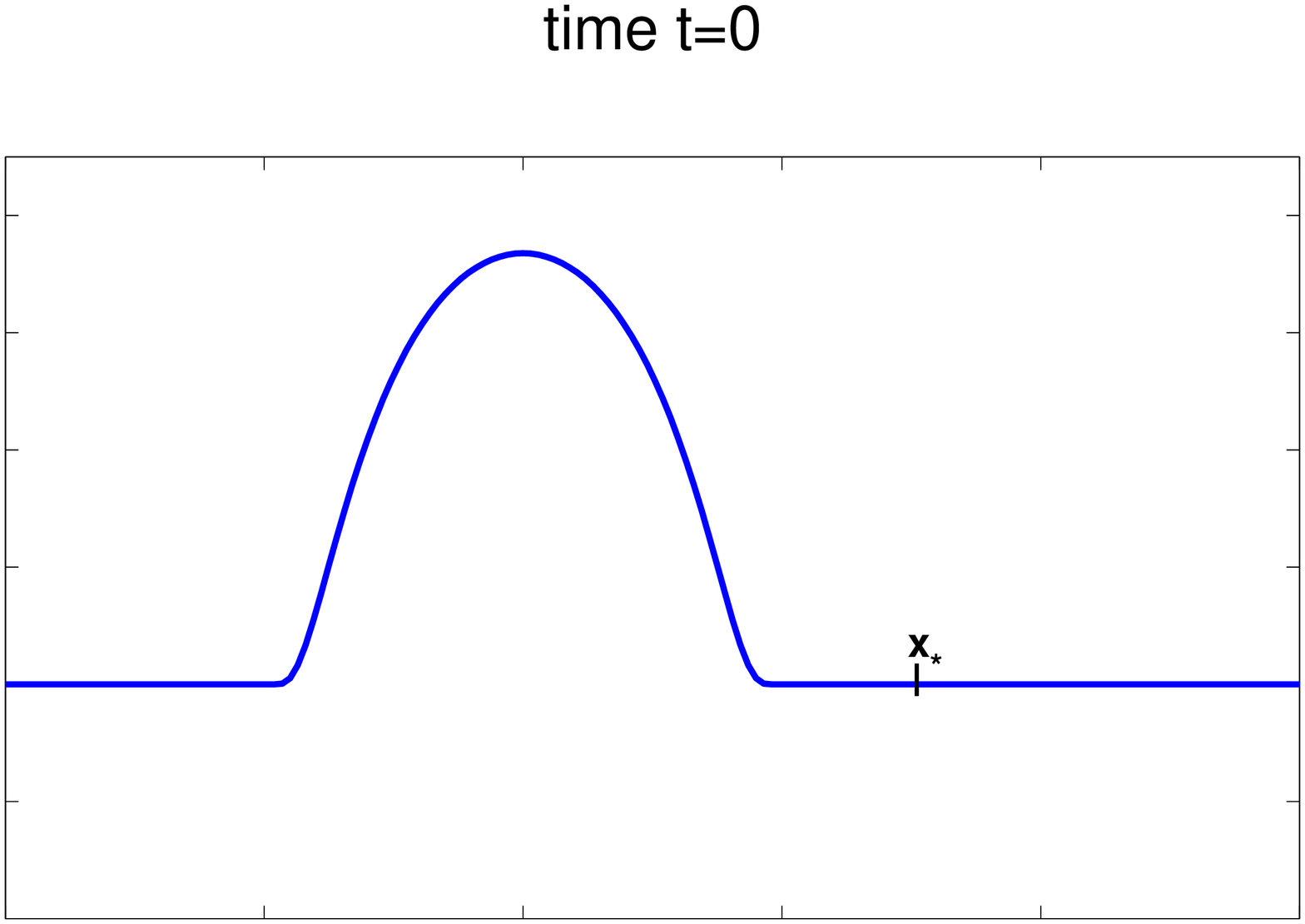} \includegraphics[scale=0.2]{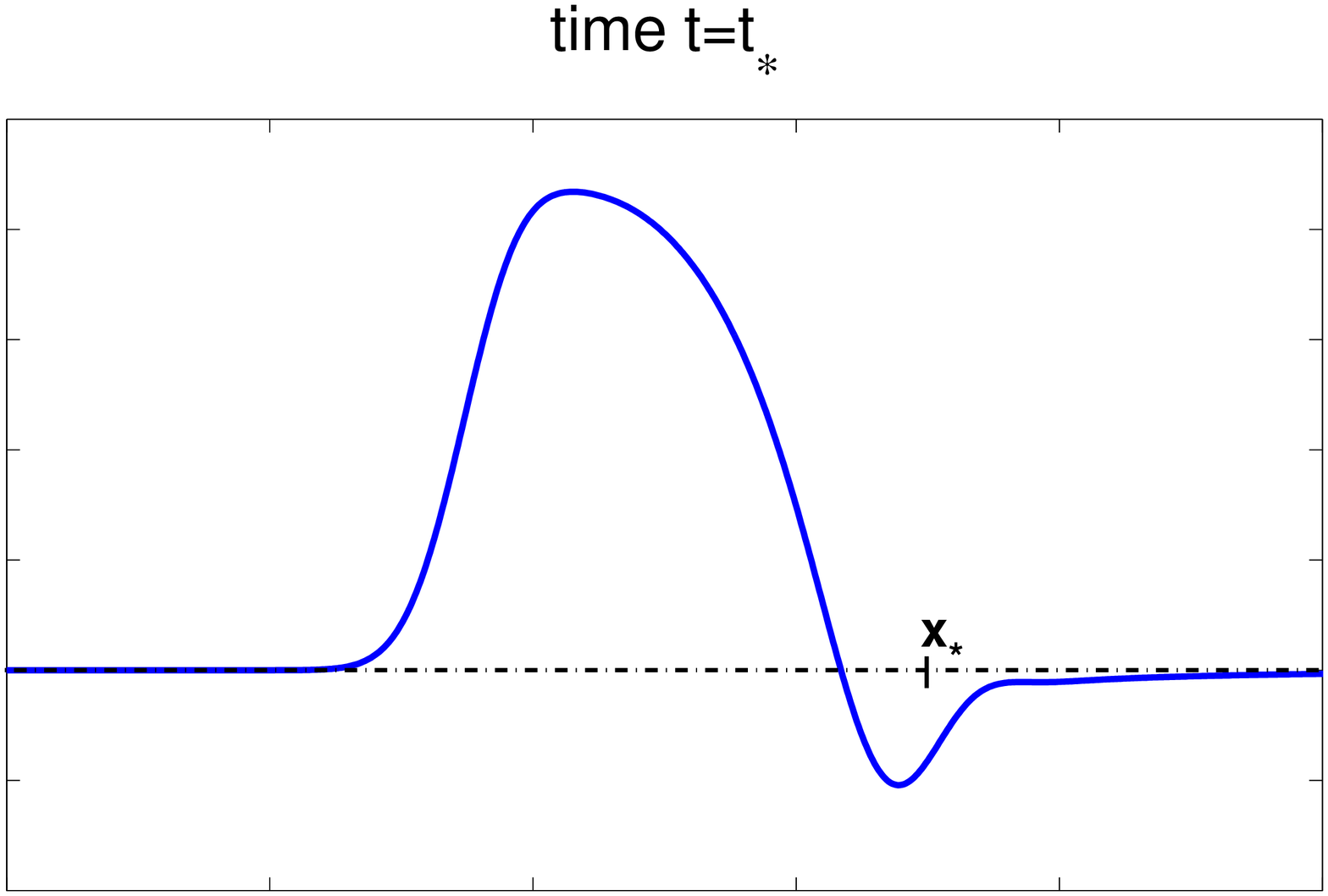}
\caption{Evolution of a dune, at $t=0$ and $t=t_\ast$. We can
observe that $u(t_\ast,x_\ast)<0$ and that $\int u(t,x)\,dx$ remains constant.} \label{non_max_princ}
\end{center}
\end{figure}
\end{remark}

\bigskip

\textbf{Notations:} In the following, we let $\F$ denote the
Fourier transform defined for $f \in L^1(\R)$ by: for all $\xi \in
\R$,
$$
\F f(\xi):= \int_\R e^{-2i \: \pi x \xi} f(x) dx.
$$
We also let $\F$ define the extension of the preceding operator
from $L^2$ to $L^2$. In the sequel, we only consider Fourier
transform with respect to (w.r.t. for short) the space variable;
in order to simplify the presentation, for any $u \in
C([0,T];L^2(\R))$, we let $\F u \in C([0,T];L^2(\R,\mathbb{C}))$
denote the function$$ t \in [0,T] \rightarrow \F(u(t,\cdot)) \in
L^2(\R,\mathbb{C}).
$$
\section{Preliminaries}\label{sect prelim}

In Subsection \ref{subsect formule}, we give the integral and the
pseudo-differential formula for $\LL$ and in Subsection
\ref{subsect noyau} we give the properties on the kernel of
$\LL-\partial_{xx}^2$.

\subsection{Integral formula for $\LL$}
\label{subsect formule}

\begin{Prop}
For all $\varphi \in \SR$ and all $x \in \R$,
\begin{equation}
\LL [\varphi](x)= C_{\LL} \int_{-\infty}^{0}
\frac{\varphi(x+z)-\varphi(x)-\varphi'(x) z}{|z|^{7/3}} \:
dz,\label{intformula}
\end{equation} \label{propintegralformula}
with $C_\LL=\frac{4}{9}$.
\end{Prop}

\begin{proof}
The proof is an easy consequence of Taylor-Poisson's formula and Fubini's
Theorem; notice that the regularity of $\varphi$ ensures the
validity of the computations that follow. We have:
\begin{eqnarray*}
\int_{-\infty}^{0} \frac{\varphi(x+z)-\varphi(x)-\varphi'(x)
z}{|z|^{7/3}} \: dz & = & \int_{-\infty}^{0} |z|^{-\frac{7}{3}}
\left(\int_0^1 (1-\tau) \varphi''(x+\tau z) z^2 d\tau \right) dz,\\
& = & \int_0^1 (1-\tau) \left(\int_{-\infty}^{0}
|z|^{-\frac{1}{3}}
\varphi''(x+\tau z)  dz \right) d\tau,\\
& = & \int_0^1 (1-\tau)\tau^{-\frac{2}{3}}
\left(\int_{0}^{+\infty} |\zeta|^{-\frac{1}{3}} \varphi''(x-\zeta)
d\zeta \right) d\tau,
\end{eqnarray*}
thanks to the change of variable $\tau z=-\zeta$. Then,
\begin{eqnarray*}
\int_{-\infty}^{0} \frac{\varphi(x+z)-\varphi(x)-\varphi'(x)
z}{|z|^{7/3}} \: dz = \int_0^1 (1-\tau) \tau^{-\frac{2}{3}}d\tau
\: \LL[\varphi](x)=\frac{9}{4} \LL[\varphi](x).
\end{eqnarray*}
The proof is now complete.
\end{proof}

\begin{Coro}
There are positive constants $a_{\LL}$ and $b_{\LL}$ such that
for all $\varphi \in \SR$ and all $\xi \in \R $,
\begin{equation}\label{formule pseudo}
\F \left( \LL [\varphi]-\varphi''\right) (\xi) = \psi_{\mathcal{I}}(\xi)  \F
\varphi (\xi),
\end{equation}
where $\psi_{\mathcal{I}}(\xi)=4 \pi^2 \xi^2-a_{\LL}
|\xi|^{\frac{4}{3}} + i \: b_{\LL} \xi |\xi|^{\frac{1}{3}}$.
\label{proppseudodiff}
\end{Coro}

\begin{proof}
We have
$$
\F \left(\LL [\varphi] \right) (\xi)=  C_{\LL}
\int_{\mathbb{R}}\int_{-\infty}^{0} e^{-2 i \pi x \xi}
\frac{\varphi(x+z)-\varphi(x)-\varphi'(x) z}{|z|^{7/3}} \: dz dx.
$$
Notice that Proposition \ref{propintegralformula} ensures that for $\varphi \in \SR$,
$\LL [\varphi] \in L^{1}(\mathbb{R})$
 and thus its Fourier
transform is well-defined. By Fubini's theorem, we can first
integrate w.r.t $x$ to deduce that
$$
\F \left(\LL [\varphi]\right) (\xi)=  C_{\LL} \int_{-\infty}^{0}
\frac {\F \left( \mathcal{T}_{-z} \varphi \right) (\xi) - \F
\varphi(\xi) - \F (\varphi')(\xi) z}{|z|^{7/3}} \: dz,
$$
where we let $\mathcal{T}_{-z} \varphi$  denote the (translated)
function $x \rightarrow  \varphi(x+z)$.
Classical formulae on Fourier transform imply that $\F \left(\LL
[\varphi]\right) (\xi) = \psi(\xi) \F \varphi (\xi) $, where
$$
\psi(\xi)= C_{\LL} \int_{-\infty}^{0} \frac {e^{2 i \pi \xi z} - 1
- 2 i \pi \xi  z}{|z|^{7/3}} \: dz.
$$
Simple computations show that
$$
\psi(\xi)=C_{\LL} \int_{-\infty}^{0} \frac { \cos{(2\pi\xi
z)-1}}{|z|^{7/3}} \: dz + i \: C_{\LL} \int_{-\infty}^{0} \frac {
 \sin{(2\pi\xi z)-2 \pi \xi z}}{|z|^{7/3}} \: dz.
$$
It is immediate that the real part of $\psi(\xi)$ is even,
non-positive, non-identically equal to $0$ and homogeneous of
degree $\frac{4}{3}$ (the last property can be seen by changing
the variable by $z'=\xi z$). Moreover, the imaginary part of
$\psi(\xi)$ is odd, negative and homogeneous of degree
$\frac{4}{3}$ on $\R_\ast^-$. There then exist positive constants
$a_{\LL}$ and $b_{\LL}$ such that
$$
\psi(\xi)=-a_{\LL} |\xi|^{\frac{4}{3}}+i \: b_{\LL} \xi |\xi|^{\frac{1}{3}}
$$
and, in particular, $\F \left(\LL
[\varphi]\right) (\xi) = \left(-a_{\LL} |\xi|^{\frac{4}{3}}+i \: b_{\LL} \xi |\xi|^{\frac{1}{3}}\right) \F \varphi (\xi) $.
Since $\F (-\varphi'')(\xi)=4 \pi^2 \xi ^2 \F \varphi(\xi)$, the proof of Corollary \ref{proppseudodiff} is complete.
\end{proof}
\begin{remark}
\begin{enumerate}
\item Since $\LL [\varphi] = \mathbf{1}_{\R_+}
|\cdot|^{-\frac{1}{3}} \ast \varphi'',$ we have $\F(\LL [\varphi])
= \F( \mathbf{1}_{\R_+} |\cdot|^{-\frac{1}{3}})\cdot
(-4\pi^2|\xi|^2)\cdot \F(\varphi)$. Elementary computations give
$\F( \mathbf{1}_{\R_+} |\cdot|^{-\frac{1}{3}}) =
\Gamma(\frac{2}{3}) \left( \frac{1}{2} -i \,\mbox{\rm{sign}}(\xi)
\frac{\sqrt{3}}{2}\right )|\xi|^{-\frac{2}{3}}$. Hence
$a_\LL=\Gamma(\frac{2}{3})\frac{1}{2}$ and
$b_\LL=\Gamma(\frac{2}{3})\frac{\sqrt{3}}{2}$. \item Let $s \in
\R$. If $\varphi \in H^s(\R),$ one can also define $\LL[\varphi]$
through its Fourier transform  by
$$\F(\LL[\varphi])(\xi) := -4\pi^2\Gamma(\frac{2}{3}) \left( \frac{1}{2} -i \,\mbox{\rm{sign}}(\xi) \frac{\sqrt{3}}{2}\right )|\xi|^{\frac{4}{3}}\cdot \F(\varphi)
$$
Thus,if $\varphi \in H^s$, we  have that $\LL[\varphi] \in
H^{s-\frac{4}{3}}$ and $||\LL[\varphi]||_{H^{s-\frac{4}{3}}} \leq
4\pi^2\Gamma(\frac{2}{3})||\varphi||_{H^{s}}.$ This implies in
particular that $\LL: H^2(\R) \rightarrow C_b(\R) \cap L^2(\R)$,
since by Sobolev embedding $H^{\frac{2}{3}} \hookrightarrow
C_b(\R) \cap L^2(\R).$
 \item Corollary \ref{proppseudodiff} implies that
$\LL-\partial_{xx}^2: C^2(\R) \cap H^2(\R) \rightarrow C(\R) \cap
L^2(\R)$ with $\LL$ which satisfies both formula
\eqref{intformula} and \eqref{formule pseudo}.
\end{enumerate}
\end{remark}

\subsection{Main properties on the kernel $K$ of $\LL-\partial_{xx}^2$}
\label{subsect noyau}

\begin{figure}[!ht!]
\begin{center}
\includegraphics[scale=0.4]{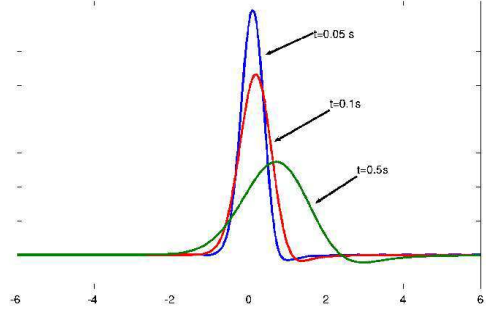}
\caption{The kernel of $\LL-\partial_{xx}^2$ for
$t=0.05, 0.1$ and $0.5$ s.} \label{kernel}
\end{center}
\end{figure}

By Corollary \ref{proppseudodiff}, we see that the semi-group
generated by $\LL-\partial^2_{xx}$ is formally given by the
convolution with the kernel (defined for $t>0$ and $x \in \R$)
$$
K(t,x)=\F^{-1} (e^{-t\psi_{\LL}})(x).
$$
\begin{Prop}
 $K(t,\cdot)$ is a $L^1$ real valued continuous
function.
\end{Prop}
\begin{proof}
 $K(t,\cdot)$ is a $L^1$ real valued continuous
function as inverse Fourier transform of a $W^{2,1}$ function with
an even real part and an odd imaginary part.
\end{proof}
In the sequel, we
only consider real valued solution of \eqref{fowlereqn}.
 We expose in Figure \ref{kernel} the evolution of $K(t,\cdot)$ for
different times. Note that $K(t,\cdot)$ is not compactly supported but that
$K(t,x) \leq  \frac{C(t)}{x^2},$ for $|x| \geq 1$ with
$C(t) = \frac{1}{4\pi^2}||\partial_{\xi\xi}^2 \F(K(t,.))(\xi) ||_{L^1}
$
\begin{Prop}
 The kernel $K$ has a non-zero negative part.
\end{Prop}
\begin{proof}
 Let us  assume  that $K$ is
nonnegative, then \label{page ref}
\begin{multline*}
|e^{-t \psi_\LL(\xi)}|  \leq ||\mathcal{F}^{-1} (e^{-t
\psi_\LL})||_{L^1(\R)}  =  \int_{\R} |K(t,.)| \\ = \int_{\R} K(t,.) =
\mathcal{F}\left( \mathcal{F}^{-1} (e^{-t \psi_\LL}) \right)(0)  =
e^{-t \psi_\LL(0)} =1
\end{multline*}
for all $\xi \in \R$; hence, since $|e^{-t
\psi_\LL(\xi)}|=e^{-t(4\pi^2|\xi|^2-a_{\LL}|\xi|^{\frac{4}{3}})}>1$
for $0<|\xi|<\frac{a_{\LL}^{\frac{3}{2}}}{8 \pi^3}$, this gives us
a contradiction.
\end{proof}
 The main consequence of this is the failure of the
maximum principle for the equation
\begin{equation}\label{equation sans burgers}
\partial_t u+\LL[u]-\partial^2_{xx} u=0;
\end{equation}
that is to say, there exists a non-negative initial condition
$u_0$ such that, for some $t>0$, $u(t,.):=K(t,.) \ast u_0$ has a
non-zero negative part, see section \ref{sect princ max} below.
 Nevertheless, $K$ enjoys many properties
similar than those one satisfied by the kernel of the heat
equation and that ensure that Equation \eqref{equation sans
burgers} has a regularizing effect on the initial condition: if
$u_0 \in L^p(\R)$ for some $p \in [1,+\infty)$, then $u$ is
$C^\infty$ for positive times, see section \ref{sect existence}.

\bigskip

Let us precise here the properties that will be needed in this
paper. Since $K(t,\cdot) \in L^1(\R)$, the family of bounded
linear operators $\left\{ u_0 \in L^2(\R) \rightarrow K(t,\cdot)
\ast u_0 \in L^2(\R) \right\}_{t>0}$ is well-defined. Moreover, it
is a strongly continuous semi-group of convolution, that is to
say:
\begin{equation}\label{semi group convolution}
\begin{array}{c}
\forall t,s>0, \: K(s,\cdot) \ast K(t,\cdot) =K(s+t,\cdot), \\
\forall u_0 \in L^2(\R), \: \lim_{t \rightarrow 0} K(t,\cdot)
\ast u_0 = u_0 \mbox{ in $L^2(\R)$.}
\end{array}
\end{equation}
Next, the kernel $K$ is smooth on $(0,+\infty) \times \R$ and we have:
\begin{equation}\label{esti gradient kernel}
\forall T>0, \: \exists \mathcal{K}_0 \mbox{ s.t. } \forall t \in
(0,T], \: ||\partial_x K(t,.)||_{L^2(\R)} \leq \mathcal{K}_0
t^{-\frac{3}{4}},
\end{equation}
\begin{equation}\label{esti gradient kernel L un}
\forall T>0, \: \exists \mathcal{K}_1 \mbox{ s.t. } \forall t \in
(0,T], \: ||\partial_x K(t,.)||_{L^1(\R)} \leq \mathcal{K}_1
t^{-\frac{1}{2}},
\end{equation}
\begin{equation}\label{semi group gradient}
\forall t,s>0, \:  K(s,\cdot) \ast \partial_x K(t,\cdot) =\partial_x
K(s+t,\cdot).
\end{equation}

\begin{proof}[Proof of these properties] The semi-group property \eqref{semi group convolution}
and \eqref{semi group gradient} are immediate consequences of
Fourier formula. Let us prove the strong continuity. By
Plancherel's formula,
\begin{multline}\label{ref strong continuity}
||K(t,\cdot) \ast u_0 -u_0||^2_{L^2(\R)}=||\F (K(t,\cdot) \ast u_0) -\F u_0||^2_{L^2(\R)}\\
=||e^{-t \psi_\LL} \F u_0 -\F u_0||^2_{L^2(\R)}=
\int_\R |e^{-t \psi_\LL}-1|^2\, |\F u_0|^2.
\end{multline}
The function $|e^{-t \psi_\LL}-1|^2\, |\F u_0|^2 $ converges
pointwise to $0$ on $\R$, as $t \rightarrow 0$. Recalling that
$\min \mbox{Re} (\psi_\LL)$ is finite, we infer that $|e^{-t
\psi_\LL}-1|^2 \, |\F u_0|^2 \leq C|\F u_0|^2$ and the dominated
convergence theorem implies that the last term of \eqref{ref
strong continuity} tends to $0$ as $t \rightarrow 0$. This
completes the proof of \eqref{semi group convolution}. Let us now
prove the estimates on the gradient. The smoothness of $K$ is an
immediate consequence of the theorem of derivation under the
integral sign applied to the definition of $K$ by Fourier
transform. We get in particular:
$$
\partial_x K(t,\cdot)=\partial_x \F^{-1}(e^{-t \psi_\LL})=\F^{-1}\left(\xi \rightarrow 2 i \pi \xi e^{-t \psi_\LL(\xi)}\right).
$$
Since the function $\xi \rightarrow 2 i \pi \xi e^{-t \psi_\LL(\xi)}$ is $L^2$, $\partial_x K(t,\cdot)$ is $L^2$ and we have:
\begin{equation*}
||\partial_x K(t,.)||^2_{L^2(\R)}=\int_\R 4 \pi^2 \xi^2 |e^{-t \psi_\LL(\xi)}|^2 d\xi
= \int_\R 4 \pi^2 \xi^2 e^{-2 t(4\pi^2|\xi|^2-a_{\LL}|\xi|^{\frac{4}{3}})} d\xi.
\end{equation*}
Let us change the variable by $\xi'=t^{\frac{1}{2}} \xi$. We get:
\begin{eqnarray*}
||\partial_x K(t,.)||^2_{L^2(\R)}
& = & t^{-\frac{3}{2}} \int_\R 4 \pi^2 |\xi'|^2 e^{-2 (4\pi^2|\xi'|^2-t^{\frac{1}{3}} a_{\LL}|\xi'|^{\frac{4}{3}})} d\xi',\\
& \leq & t^{-\frac{3}{2}} \int_\R 4 \pi^2 |\xi'|^2 e^{-2 (4\pi^2|\xi'|^2-T^{\frac{1}{3}} a_{\LL}|\xi'|^{\frac{4}{3}})} d\xi',
\end{eqnarray*}
for all $t \in (0,T]$. The proof of \eqref{esti gradient kernel}
is now complete. To prove \eqref{esti gradient kernel L un}, we
have to derive a ''homogeneity-like'' property for $K$. Easy
computations show that
\begin{eqnarray*}
K(t,x) & = & \int_\R e^{2 i \: \pi x \xi} e^{-t \psi_\LL(\xi)} d\xi,\\
& = & \int_\R e^{2 i \: \pi x \xi} e^{-t (4 \pi^2 |\xi|^2-a_\LL
|\xi|^\frac{4}{3}+i \: b_\LL \xi |\xi|^\frac{1}{3})} d\xi,\\
& = & t^{-\frac{1}{2}} \int_\R e^{2 i \: \pi (t^{-\frac{1}{2}}  x)
\xi'} e^{- (4 \pi^2 |\xi'|^2-t^{\frac{1}{3}}a_\LL
|\xi'|^\frac{4}{3}+i \:t^{\frac{1}{3}} b_\LL \xi'
|\xi'|^\frac{1}{3})} d\xi',
\end{eqnarray*}
by changing the variable by $\xi'=t^{\frac{1}{2}} \xi$. Then,
\begin{eqnarray*}
K(t,x) & = & t^{-\frac{1}{2}} \int_\R e^{2 i \: \pi
(t^{-\frac{1}{2}}  x) \xi'} e^{- (4 \pi^2 |\xi'|^2-a_\LL
|\xi'|^\frac{4}{3}+i \: b_\LL \xi' |\xi'|^\frac{1}{3})}
e^{-(1-t^{\frac{1}{3}})(a_\LL|\xi'|^\frac{4}{3}-i \: b_\LL
\xi'|\xi'|^\frac{1}{3})} d\xi',\\
 & = & t^{-\frac{1}{2}}\int_\R e^{2 i \: \pi
(t^{-\frac{1}{2}}  x) \xi'} e^{- \psi_\LL(\xi')}
e^{-(1-t^{\frac{1}{3}})(a_\LL|\xi'|^\frac{4}{3}-i \: b_\LL
\xi'|\xi'|^\frac{1}{3})} d\xi'.
\end{eqnarray*}
For $t<1$, define $G((1-t^{\frac{1}{3}}),\cdot):=
\F^{-1}\left(e^{-(1-t^{\frac{1}{3}})(a_\LL|\xi'|^\frac{4}{3}-i \: b_\LL
\xi'|\xi'|^\frac{1}{3})} \right)$. It is readily seen that $G$ is
$L^1$ as inverse Fourier transform of a $W^{2,1}$ function.
Moreover, for $t_0 \in (0,1)$ and all $t \in (0,t_0]$,
$$
||G((1-t^{\frac{1}{3}}),\cdot)||_{L^1(\R)} \leq C
\left|\left|e^{-(1-t^{\frac{1}{3}})(a_\LL|\cdot|^\frac{4}{3}-i \: b_\LL
\cdot|\cdot|^\frac{1}{3})}\right| \right|_{W^{2,1}(\R,\C)} \leq
C(t_0),
$$
where $C(t_0)$ only depends on $t_0$. 
Classical
formula on Fourier transform then give:
$$
K(t,x) = t^{-\frac{1}{2}} \left(K(1,\cdot) \ast
G((1-t^\frac{1}{3}),\cdot) \right)(t^{-\frac{1}{2}}x).
$$
Observe now that $\partial_x K(1,\cdot)=\F^{-1} \left(\xi
\rightarrow 2i \: \xi \pi e^{-\psi_\LL(\xi)} \right)$ is $L^1$ as
inverse Fourier transform of a $W^{2,1}$ function. Then,
$$
\partial_x K(t,x) = t^{-1} \left(\partial_x K(1,\cdot) \ast
G((1-t^\frac{1}{3}),\cdot) \right)(t^{-\frac{1}{2}}x)
$$
is $L^1$ and its $L^1$ norm can be computed by the change of
variable $x'=t^{-\frac{1}{2}} x$ as follows:
$$
||\partial_x K(t,\cdot)||_{L^1(\R)}= t^{-\frac{1}{2}} ||
\partial_x K(1,\cdot) \ast
G((1-t^\frac{1}{3}),\cdot) ||_{L^1(\R)} \leq t^{-\frac{1}{2}} ||
\partial_x K(1,\cdot) ||_{L^1(\R)} C(t_0),
$$
for any $t \in (0,t_0]$. Since
$$
||\partial_x K(t,\cdot)||_{L^1(\R)} \leq C||\xi \rightarrow 2i \:
\xi \pi e^{-t \psi_\LL(\xi)} ||_{W^{2,1}(\R,\C)} \leq C(t_0,T),
$$
for all $t \in [t_0,T]$, the proof of \eqref{esti gradient kernel
L un} is now complete.
\end{proof}

\begin{remark}\label{omega 0}
For any $u_0 \in L^2(\R)$ and $t>0$,
\begin{equation}\label{esti omega 0}
||K(t,\cdot) \ast u_0||_{L^2(\R)} \leq e^{\omega_0 t}
||u_0||_{L^2(\R)},
\end{equation}
where $\omega_0=-\min \mbox{\emph{Re}} (\psi_\LL)$.
\end{remark}
\begin{proof}
 This is readily established with Plancherel's formula, like in \eqref{ref strong continuity}.
\end{proof}

\section{Duhamel's formula}
\label{sect duhamel}

Using Fourier transform and Corollary \ref{proppseudodiff}, we
formally see that any solution to \eqref{fowlereqn} satisfies
Duhamel's formula \eqref{formule duhamel} (see also the proof of
Lemma \ref{prop: reg bord}, which justifies the computations).
This observation is the starting point of the definition of mild
solution below.
\begin{definition}
Let $T>0$ and $u_0 \in L^2(\R)$. We say that $u \in
L^\infty((0,T);L^2(\R))$ is a mild solution to \eqref{fowlereqn}
if for a.e. $t \in (0,T)$,
\begin{equation}\label{formule duhamel}
u(t,\cdot)=K(t,\cdot) \ast u_0 - \frac{1}{2} \int_{0}^{t}
\partial_x K (t-s,\cdot) \ast u^2(s,\cdot) ds.
\end{equation}
\label{def Duhamel}
\end{definition}

The following proposition shows that all the terms in
\eqref{formule duhamel} are well-defined and that Equation \eqref{fowlereqn} generates a (non-linear) semi-group.
\begin{Prop}\label{prop: cont semi groupe}
Let $T>0$, $u_0 \in L^2(\R)$ and $v \in L^\infty((0,T);L^1(\R))$. Then, the function
\begin{equation}\label{def: theta}
u: t \in (0,T] \rightarrow K(t,\cdot) \ast u_0 - \frac{1}{2} \int_{0}^{t}
\partial_x K (t-s,\cdot) \ast v(s,\cdot)ds \in L^2(\R),
\end{equation}
is well-defined and belongs to $C([0,T];L^2(\R))$ (being
extended at $t=0$ by the value $u(0,.)=u_0$).\\
(Semi-group property) Moreover, for all $t_0 \in (0,T)$ and all $t \in [0,T-t_0]$,
\begin{equation*}
u(t_0+t,.)=K(t,\cdot) \ast u(t_0,\cdot) - \frac{1}{2} \int_{0}^{t}
\partial_x K (t-s,\cdot) \ast v(t_0+s,\cdot) ds.
\end{equation*}
\end{Prop}
\begin{proof}
By \eqref{semi group convolution}, it is  classical that the
function $t \in (0,T] \rightarrow K(t,\cdot) \ast u_0 \in L^2(\R)$
is continuous and can be continuously extended by the value
$u(0,\cdot)=u_0$ at $t=0$. What is left to prove  is thus the continuity of the function
$$
w: t \in [0,T] \rightarrow \int_0^t \partial_x K (t-s,\cdot) \ast v(s,\cdot) ds \in L^2(\R).
$$
Let us extend $\partial_x K$ and $v$ for all times the following way:
$$
\mathcal{H}(t,\cdot):= \left\{
\begin{array}{cc}
\partial_x K(t,\cdot) & \mbox{ if $t>0$,}\\
0 & \mbox{ if not}
\end{array}
\right. \quad \mbox{ and } \quad \mathcal{V}(t,\cdot):= \left\{
\begin{array}{cc}
v(t,\cdot) & \mbox{ if $t \in (0,T)$,}\\
0 & \mbox{ if not.}
\end{array}
\right.
$$
Then we have
$$
w(t,\cdot)=\int_\R \mathcal{H}(t-s,\cdot) \ast \mathcal{V}(s,\cdot)ds.
$$
It is immediate that $\mathcal{V} \in L^\infty(\R;L^1(\R))$. Moreover, \eqref{esti gradient kernel} implies that
\begin{equation}\label{esti tech 1}
|| \mathcal{H}(t,\cdot)||_{L^2(\R)} \leq \mathbf{1}_{\{0<t<T\}}
\mathcal{K}_0 t^{-\frac{3}{4}}
\end{equation}
and it follows that $\mathcal{H} \in L^1(\R;L^2(\R))$.
Young's Inequalities imply that for all $t \in \R$
\begin{eqnarray}
\int_\R ||\mathcal{H}(t-s,\cdot) \ast \mathcal{V}(s,\cdot)||_{L^2(\R)} ds & \leq &
\int_\R ||\mathcal{H}(t-s,\cdot)||_{L^2(\R)} ||\mathcal{V}(s,\cdot)||_{L^1(\R)} ds, \nonumber \\
& \leq & ||\mathcal{H}||_{L^1(\R;L^2(\R))} \: ||\mathcal{V}||_{L^\infty(\R;L^1(\R))}. \label{esti tech 2}
\end{eqnarray}
This implies, in particular, that the function $w$ is well-defined. Let us now take $t,s \in \R$ and define
$$
I:=\left| \left|\int_\R \mathcal{H}(t-\tau,\cdot) \ast
\mathcal{V}(\tau) d\tau-\int_\R \mathcal{H}(s-\tau,\cdot) \ast
\mathcal{V}(\tau) d\tau \right| \right|_{L^2(\R)}.
$$
We have
\begin{eqnarray*}
I & \leq & \int_\R \left| \left| \left(\mathcal{H}(t-\tau,\cdot)-\mathcal{H}(s-\tau,\cdot) \right) \ast \mathcal{V}(\tau) \right| \right|_{L^2(\R)} d\tau,\\
& \leq & \int_\R ||\mathcal{H}(t-\tau,\cdot)-\mathcal{H}(s-\tau,\cdot)||_{L^2(\R)} ||\mathcal{V}(\tau)||_{L^1(\R)}d\tau,
\end{eqnarray*}
thanks to Young's Inequalities. It follows that
$$
I \leq \int_\R ||\mathcal{H}(t-\tau,\cdot)-\mathcal{H}(s-\tau,\cdot)||_{L^2(\R)} d\tau \: ||\mathcal{V}||_{L^\infty(\R;L^1(\R))}.
$$
Since the translation are continuous in $L^1(\R;L^2(\R))$, we see
that $I \rightarrow 0$ as $|t-s| \rightarrow 0$. In particular,
the function $w$ is continuous and this completes the proof of the
continuity of $u$.

Let us now prove the semi-group property. By \eqref{semi group
convolution} and \eqref{semi group gradient}, we infer that
\begin{eqnarray*}
u(t_0+t,\cdot) & = & K(t,\cdot) \ast K(t_0,\cdot) \ast u_0-\frac{1}{2} \int_0^{t_0} \partial_x K(t+t_0-s,\cdot) \ast v(s,\cdot) ds\\
 & & \quad \quad -\frac{1}{2}  \int_{t_0}^{t+t_0} \partial_x K(t+t_0-s,\cdot) \ast v(s,\cdot) ds,\\
 & = & K(t_0,\cdot) \ast K(t,\cdot) \ast u_0-\frac{1}{2} \int_0^{t_0} K(t,\cdot) \ast \partial_x K(t_0-s,\cdot) \ast v(s,\cdot) ds\\
 & & \quad \quad -\frac{1}{2} \int_{0}^{t} \partial_x K(t-s',\cdot) \ast v(t_0+s',\cdot) ds',
\end{eqnarray*}
thanks to the change of variable $s'=s-t_0$ to compute the last integral term. Then,
\begin{eqnarray*}
u(t_0+t,\cdot) & = & K(t,\cdot) \ast K(t_0,\cdot) \ast u_0-K(t,\cdot) \ast \frac{1}{2} \int_0^{t_0} \partial_x K(t_0-s,\cdot) \ast v(s,\cdot) ds\\
 & & \quad \quad -\frac{1}{2}  \int_{0}^{t} \partial_x K(t-s',\cdot) \ast v(t_0+s',\cdot) ds',\\
 & = & K(t,\cdot) \ast \left(K(t_0,\cdot) \ast u_0-\frac{1}{2} \int_0^{t_0} \partial_x K(t_0-s,\cdot) \ast v(s,\cdot) ds \right)\\
& & \quad \quad  -\frac{1}{2}  \int_{0}^{t} \partial_x K(t-s',\cdot) \ast v(t_0+s',\cdot) ds',\\
 & = & K(t,\cdot) \ast u(t_0,\cdot)-\frac{1}{2} \int_{0}^{t} \partial_x K(t-s',\cdot) \ast v(t_0+s',\cdot) ds'.
\end{eqnarray*}
The proof of the semi group property is now complete.
\end{proof}
\begin{remark}\label{rem esti L 2}
For $ v \in
L^\infty((0,T);L^1(\R))$, $u \in C([0,T];L^2(\R))$ defined in \eqref{def: theta}
satisfies:
\begin{equation}\label{esti: duhamel}
||u||_{C([0,T];L^2(\R))} \leq e^{\omega_0 T}
||u_0||_{L^2(\R)}+2\mathcal{K}_0 T^{\frac{1}{4}}
||v||_{L^\infty((0,T);L^1(\R))}.
\end{equation}
\end{remark}
\begin{proof}
 Indeed, with \eqref{esti tech 1} and \eqref{esti tech 2}, we  estimate
the integral term of \eqref{def: theta} and with  \eqref{esti omega
0}, we   estimate the $L^2$ norm of $K(t,\cdot) \ast u_0$.
\end{proof}
\section{Uniqueness of a solution}
\label{sect uniqueness}
Let us state a lemma  that will be needed later.
\begin{lemme}\label{rem: contraction}
 Let
$T>0$, $u_0 \in L^2(\R)$. For
$i=1,2$, let  $v_i \in L^\infty((0,T);L^1(\R))$  and define $u_i \in C([0,T];L^2(\R))$  as in Proposition
\ref{prop: cont semi groupe} by
$$u_i(t,\cdot):= K(t,\cdot) \ast u_0 - \frac{1}{2} \int_{0}^{t}
\partial_x K (t-s,\cdot) \ast v_i(s,\cdot)ds.
$$
 Then we have the estimate
\begin{equation}\label{esti: contraction bonne}
|| u_1-u_2||_{C([0,T];L^2(\R))} \leq 2 \mathcal{K}_0
T^{\frac{1}{4}} \: ||v_1-v_2||_{L^\infty((0,T);L^1(\R))}.
\end{equation}
\end{lemme}
\begin{proof}
  For
all $t \in [0,T]$, we have
\begin{equation*}
u_1(t,\cdot)-u_2(t,\cdot)= - \frac{1}{2} \int_{0}^{t} \partial_x K
(t-s,\cdot) \ast (v_1(s,\cdot)-v_2(s,\cdot)) ds.
\end{equation*}
Hence,
\begin{eqnarray}
|| u_1(t,\cdot)-u_2(t,\cdot)||_{L^2(\R)} & =  & \frac{1}{2} \left|
\left|
 \int_{0}^{t} \partial_x K (t-s,\cdot)
\ast (v_1(s,\cdot)-v_2(s,\cdot)) ds \right| \right|_{L^2(\R)} , \nonumber\\
& \leq & \frac{1}{2} \int_{0}^{t} ||
\partial_x K (t-s,.) \ast (v_1(s,.)-v_2(s,.))||_{L^2(\R)}\, ds.
\label{uniq-1}
\end{eqnarray}
By (\ref{esti gradient kernel}),
\begin{eqnarray*}
|| \partial_x K (t-s,.) \ast
(v_1(s,\cdot)-v_2(s,\cdot))||_{L^2(\R)}
        & \leq & || \partial_x K(t-s,\cdot) \, ||_{L^2(\R)} || v_1(s,\cdot)-v_2(s,\cdot)||_{L^1(\R)} \\
        & \leq &   \mathcal{K}_0 (t-s)^{-\frac{3}{4}} ||
        v_1(s,\cdot)-v_2(s,\cdot)||_{L^1(\R)}.
\end{eqnarray*}
Inequality \eqref{uniq-1} then gives
\begin{eqnarray*}
|| u_1(t,\cdot)-u_2(t,\cdot)||_{L^2(\R)}  & \leq & \frac{
\mathcal{K}_0}{2} \int_{0}^{t} (t-s)^{-\frac{3}{4}} ds \:
||v_1-v_2||_{L^\infty((0,t);L^1(\R))},\\
& = & 2 \mathcal{K}_0 t^{\frac{1}{4}} \:
||v_1-v_2||_{L^\infty((0,t);L^1(\R))}.
\end{eqnarray*}
In particular, for all $s \in [0,t]$
$$
|| u_1(s,\cdot)-u_2(s,\cdot)||_{L^2(\R)} \leq 2\mathcal{K}_0
s^{\frac{1}{4}} \: ||v_1-v_2||_{L^\infty((0,s);L^1(\R))} \leq 2
\mathcal{K}_0 t^{\frac{1}{4}} \:
||v_1-v_2||_{L^\infty((0,t);L^1(\R))}
$$
and we have proved that
\begin{equation}\label{esti tech contraction}
||u_1-u_2||_{C([0,t];L^2(\R))} \leq 2 \mathcal{K}_0
t^{\frac{1}{4}} \: ||v_1-v_2||_{L^\infty((0,t);L^1(\R))}.
\end{equation}
\end{proof}
\begin{Prop}\label{prop: uniqueness}
Let $T>0$ and $u_0 \in L^2(\R)$. There exists at most one $u \in
L^\infty((0,T);L^2(\R))$ which is a mild solution to
\eqref{fowlereqn}.
\end{Prop}

\begin{proof}
Let $u,v \in L^\infty((0,T);L^2(\R))$ be two mild solutions.
Let $t \in [0,T]$.
With Lemma \ref{rem: contraction} applied to $v_1= u^2$ and $v_2=v^2$, we get
\begin{equation}\label{esti tech contraction}
||u-v||_{C([0,t];L^2(\R))} \leq 2 \mathcal{K}_0 t^{\frac{1}{4}} \:
||u^2-v^2||_{L^\infty((0,t);L^1(\R))}.
\end{equation}
Since $||u^2-v^2||_{L^\infty((0,t);L^1(\R))} \leq M
||u-v||_{C([0,t];L^2(\R))}$ with
$M=||u||_{C([0,T],L^2(\R))}+||v||_{C([0,T],L^2(\R))}$, we get:
\begin{equation*}
||u-v||_{C([0,t];L^2(\R))} \leq 2M \mathcal{K}_0 t^{\frac{1}{4}}
\: ||u-v||_{C([0,t];L^2(\R))}.
\end{equation*}
We then have established that $u=v$ on $[0,t]$ for any $t \in
(0,T]$ such that $t<(2M \mathcal{K}_0)^{-4}.$  Notice that since
$u$ and $v$ are continuous with values in $L^2$, $u=v$ on
$[0,T_\ast]$ with $T_\ast=(2M \mathcal{K}_0)^{-4}>0 $. To prove
that $u=v$ on $[0,T]$, let us define $t_0:=\sup \{t \in (0,T]
\mbox{ s.t. $u=v$ on $[0,t]$} \}$ and let us assume that $t_0 \neq
T$. The continuity of $u$ and $v$ implies that
$u(t_0,\cdot)=v(t_0,\cdot)$. The semi-group property of
Proposition \ref{prop: cont semi groupe} thus implies that
$u(t_0+\cdot,\cdot)$ and $v(t_0+\cdot,\cdot)$ are mild solutions
of \eqref{fowlereqn}with the same initial condition; that is to
say $u(t_0+0,\cdot)=v(t_0+0,\cdot)$. The first step of the proof
then implies that $u(t_0+\cdot,\cdot)=v(t_0+\cdot,\cdot)$ on
$[0,\min \{T_\ast,T-t_0 \}]$; hence, we get a contradiction with
the definition of $t_0$ and we deduce that $t_0=T$. The proof of
the uniqueness is now complete.
\end{proof}

\section{Existence of a regular solution}
\label{sect existence}

This section is devoted to the proof of the existence of a
solution $u \in C^{1,2}((0,T] \times \R)$to \eqref{fowlereqn};
that is to say, $u$ is $C^2$ in space and $C^1$ in time. We first
need the following technical result:
\begin{lemme}\label{lemme reg}
Let $u_0 \in L^2(\R)$ and $T>0$. Let $v \in C([0,T];L^1(\R)) \cap
C((0,T];W^{1,1}(\R))$ that satisfies
\begin{equation}\label{hyp v reg}
\sup_{t \in (0,T]} t^\frac{1}{2} ||\partial_{x}
v(t,\cdot)||_{L^1(\R)}<+\infty.
\end{equation}
Let $u \in C([0,T];L^2(\R))$ be the function defined in
\eqref{def: theta}. Then, $u \in C((0,T];H^1(\R))$ with
\begin{equation}\label{esti L2 gradient duhamel}
\sup_{t \in (0,T]} t^{\frac{1}{2}} ||\partial_x
u(t,\cdot)||_{L^2(\R)} \leq \mathcal{K}_1 ||u_0||_{L^2(\R)} +
\frac{\mathcal{K}_0 I}{2} \; {T^\frac{1}{4}} \sup_{t \in (0,T]}
t^\frac{1}{2} ||\partial_{x} v(t,\cdot)||_{L^1(\R)},
\end{equation}
where $I$ is a constant equal to $\int_0^1 (1-s)^{-\frac{3}{4}}
s^{-\frac{1}{2}} ds = B(1/2,1/4)$, $B$ being the beta function.\\
Moreover, let $v_i \in C([0,T];L^1(\R)) \cap C((0,T];W^{1,1}(\R))$ satisfy \eqref{hyp v
reg} and define $u_i$ by \eqref{def: theta} (with $u$ and $v$
replaced, respectively, by $u_i$ and $v_i$)  for $i=1,2$. Then,
\begin{equation}\label{contrac reg}
\sup_{t \in (0,T]} t^{\frac{1}{2}}||\partial_x (u_1-
u_2)(t,\cdot)||_{L^2(\R)} \leq \frac{\mathcal{K}_0 I}{2} \;
T^{\frac{1}{4}} \sup_{t \in (0,T]} t^\frac{1}{2} ||\partial_{x}
(v_1-v_2)(t,\cdot)||_{L^1(\R)}.
\end{equation}
\end{lemme}
\begin{proof}
Recall that Proposition \ref{prop: cont semi groupe} ensures that
$u \in C([0,T];L^2(\R))$. It is easy to check that the 
distribution derivative of $u$ w.r.t. the space variable
satisfies: for any $t \in (0,T]$,
\begin{equation*}
\partial_x u(t,\cdot)= \partial_x K(t,\cdot) \ast u_0 -
\frac{1}{2}\int_{0}^{t}
\partial_x K (t-s,\cdot) \ast \partial_x v(s,\cdot) ds.
\end{equation*}
Let us verify that all the terms are well-defined in $L^2$. Since
$\partial_x K(t,\cdot) \in L^1(\R)$, it is obvious that
$\partial_x K(t,\cdot) \ast u_0 \in L^2(\R)$. Moreover, define
$$
w(t,\cdot):=\frac{1}{2}\int_{0}^{t}
\partial_x K (t-s,\cdot) \ast \partial_x v(s,\cdot) ds.
$$
Young's Inequalities and \eqref{esti gradient kernel} give
\begin{eqnarray}
||\partial_x K (t-s,\cdot) \ast \partial_x v(s,\cdot)||_{L^2(\R)}
& \leq &
||\partial_x K(t-s,\cdot)||_{L^2(\R)} ||\partial_x v(s,\cdot)||_{L^1(\R)}, \nonumber \\
& = & ||\partial_x K(t-s,\cdot)||_{L^2(\R)} s^{-\frac{1}{2}}
s^{\frac{1}{2}} ||\partial_x v(s,\cdot)||_{L^1(\R)},\nonumber  \\
& \leq & \mathcal{K}_0(t-s)^{-\frac{3}{4}} s^{-\frac{1}{2}}
\sup_{\tau \in (0,T]} \tau^\frac{1}{2} ||\partial_{x}
v(\tau,\cdot)||_{L^1(\R)}. \label{tech exists 66}
\end{eqnarray}
Since  $\int_0^t (t-s)^{-\frac{3}{4}} s^{-\frac{1}{2}}\, ds < \infty$,
by \eqref{hyp v reg} we deduce that $w(t,\cdot)$ is well-defined
in $L^2$ and thus for all $t \in (0,T]$, $\partial_x u(t,\cdot) \in L^2(\R)$. Let us now prove that $\partial_x u$ is continuous on
$(0,T]$ with values in $L^2$. For $\delta
>0$ and $t \in (0,T]$, define
$$
w_\delta(t,\cdot):=\frac{1}{2}\int_{0}^{t}
\partial_x K (t-s,\cdot) \ast \left(\mathbf{1}_{\{s> \delta\}} \partial_x v(s,\cdot)\right) ds.
$$
Since $\mathbf{1}_{\{s> \delta\}} \partial_x v(s,\cdot) \in L^\infty ( [0,T];L^1(\R))$,
Proposition \ref{prop: cont semi groupe} ensures that $w_\delta$
is continuous on $[0,T]$ with values in $L^2$. Moreover, for any
$t_0 \in (0,T]$, $\delta \leq t_0$ and $t \in [t_0,T]$,
\begin{eqnarray*}
||w(t,\cdot)-w_\delta(t,\cdot)||_{L^2(\R)} & \leq &
\frac{1}{2}\int_{0}^{\delta}
||\partial_x K (t-s,\cdot) \ast \partial_x v(s,\cdot)||_{L^2(\R)} ds,\\
 & \leq &  \frac{\mathcal{K}_0}{2}
\int_{0}^{\delta} (t-s)^{-\frac{3}{4}} s^{-\frac{1}{2}} ds \sup_{s
\in (0,T]} s^\frac{1}{2} ||\partial_{x} v(s,\cdot)||_{L^1(\R)}
\quad
\mbox{by \eqref{tech exists 66},}\\
& \leq &  \frac{\mathcal{K}_0}{2}  \int_{0}^{\delta}
(t_0-s)^{-\frac{3}{4}} s^{-\frac{1}{2}} ds \sup_{s \in (0,T]}
s^\frac{1}{2} ||\partial_{x} v(s,\cdot)||_{L^1(\R)}.
\end{eqnarray*}
It follows that
$$
\sup_{t \in [t_0,T]} ||w(t,\cdot)-w_\delta(t,\cdot)||_{L^2(\R)}
\leq  \frac{\mathcal{K}_0}{2}  \int_{0}^{\delta}
(t_0-s)^{-\frac{3}{4}} s^{-\frac{1}{2}} ds \sup_{s \in (0,T]}
s^\frac{1}{2} ||\partial_{x} v(s,\cdot)||_{L^1(\R)} \rightarrow 0,
$$
as $\delta \rightarrow 0$. We deduce that $w \in C((0,T];L^2(\R))$
as local uniform limit of continuous functions. Moreover,
$$
\partial_x K(t,\cdot) \ast u_0 =\F^{-1} \left(\xi \rightarrow 2i \: \pi \xi e^{-t \psi_\LL(\xi)} \F u_0(\xi)
\right).
$$
The dominated convergence theorem immediately implies that for any
$t_0>0$,
$$
\int_\R  4 \pi^2 |\xi|^2 \left|e^{-t \psi_\LL(\xi)}-e^{-t_0
\psi_\LL(\xi)}\right|^2 | \F u_0(\xi)|^2 d\xi \rightarrow 0, \quad
\mbox{as $t \rightarrow t_0$}.
$$
This means that $t >0 \rightarrow \left(\xi \rightarrow 2i \: \pi
\xi e^{-t \psi_\LL(\xi)} \F u_0\right) \in L^2(\R)$ is continuous
and, since $\F$ is an isometry of $L^2$, we deduce that $t > 0
\rightarrow
\partial_x K(t,\cdot) \ast u_0 \in L^2(\R)$ is continuous. We then have established that
$\partial_x u \in C((0,T];L^2(\R))$. Let us now estimate how the
$L^2$ norm of $\partial_x u$ can explode at $t=0$. By
\eqref{tech exists 66},
\begin{multline*}
||w(t,\cdot)||_{L^2(\R)} \leq  \frac{\mathcal{K}_0}{2}  \int_0^t
(t-s)^{-\frac{3}{4}} s^{-\frac{1}{2}} ds \sup_{\tau \in (0,T]}
\tau^\frac{1}{2} ||\partial_{x} v(\tau,\cdot)||_{L^1(\R)} \\ =
\frac{\mathcal{K}_0 I}{2} t^{-\frac{1}{4}} \sup_{\tau \in (0,T]}
\tau^\frac{1}{2} ||\partial_{x} v(\tau,\cdot)||_{L^1(\R)},
\end{multline*}
where $I=\int_0^1 (1-s')^{-\frac{3}{4}} s'^{-\frac{1}{2}} ds' = B(1/2, 1/4)$;
notice that the last integral term has been computed with the help
of the change of variable $s'=\frac{s}{t}$. Moreover, \eqref{esti
gradient kernel L un} and Young's Inequalities imply that
$$
||\partial_x K(t,\cdot) \ast u_0||_{L^2(\R)} \leq \mathcal{K}_1
t^{-\frac{1}{2}}  ||u_0||_{L^2(\R)} .
$$
We deduce that for any $t \in (0,T]$,
\begin{equation*}
||\partial_x u(t,\cdot)||_{L^2(\R)} \leq \mathcal{K}_1
t^{-\frac{1}{2}} ||u_0||_{L^2(\R)} + \frac{\mathcal{K}_0 I}{2}
t^{-\frac{1}{4}} \sup_{s \in (0,T]} s^\frac{1}{2} ||\partial_{x}
v(s,\cdot)||_{L^1(\R)},
\end{equation*}
which implies immediately \eqref{esti L2 gradient duhamel}.

Let us now prove \eqref{contrac reg}. For any $t \in (0,T]$,
\begin{eqnarray*}
||\partial_x (u_1-u_2)(t,\cdot)||_{L^2(\R)} & \leq & \frac{1}{2}
\int_{0}^{t} ||
\partial_x K (t-s,.) \ast \partial_x (v_1-v_2)(s,\cdot)||_{L^2(\R)} ds,\\
& \leq & \frac{\mathcal{K}_0}{2} \int_{0}^{t}
(t-s)^{-\frac{3}{4}}s^{-\frac{1}{2}}ds \sup_{s \in (0,T]} s^{\frac{1}{2}} ||\partial_x (v_1-v_2)(s,\cdot)||_{L^1(\R)}, \\
& = & \frac{\mathcal{K}_0 I}{2} t^{-\frac{1}{4}} \sup_{s \in
(0,T]} s^{\frac{1}{2}} ||\partial_x
(v_1-v_2)(s,\cdot)||_{L^1(\R)},
\end{eqnarray*}
which implies immediately \eqref{contrac reg}.
\end{proof}

\begin{remark}\label{rem for grad}
Let $u_0, T, v$ and $u$ that satisfy assumptions of Lemma
\ref{lemme reg}. Then, we
have established that for any $t \in (0,T]$,
\begin{equation*}
\partial_x u(t,\cdot)= \partial_x K(t,\cdot) \ast u_0 -
\frac{1}{2}\int_{0}^{t}
\partial_x K (t-s,\cdot) \ast \partial_x v(s,\cdot) ds.
\end{equation*}
\end{remark}

Let us now prove the local-in-time existence of a regular solution.
\begin{Prop}\label{prop: exist local}
Let $u_0 \in L^2(\R)$. There exists $T_\ast>0$ that only depends on
$||u_0||_{L^2(\R)}$ such that \eqref{fowlereqn} admits a (unique)
mild solution $u \in C([0,T_\ast];L^2(\R)) \cap
C((0,T_\ast];H^2(\R))$ on $(0,T_\ast)$such that
$$
\sup_{t \in (0,T_\ast]} t^\frac{1}{2} ||\partial_{x}
u(t,\cdot)||_{L^2(\R)} < +\infty \quad \mbox{and} \quad \sup_{t
\in (0,T_\ast]} t ||\partial_{xx}^2 u(t,\cdot)||_{L^2(\R)} <
+\infty.
$$
Moreover, $u$ belongs to
$C^{1,2}((0,T_\ast] \times \R)$ and satisfies the PDE in
\eqref{fowlereqn} in the classical sense.
\end{Prop}

\begin{proof}
We use a contracting fixed point theorem. For $u \in
C([0,T_\ast];L^2(\R)) \cap C((0,T_\ast];H^1(\R))$, define the norm
\begin{equation}\label{norme regulari}
|||u|||:= ||u||_{C([0,T_\ast];L^2(\R))}+\sup_{t \in (0,T_\ast]}
t^\frac{1}{2} ||\partial_{x} u(t,\cdot)||_{L^2(\R)}.
\end{equation}
Define the space
$$
X:= \left\{u \in C([0,T_\ast];L^2(\R)) \cap C((0,T_\ast];H^1(\R))
\mbox{ s.t. $u(0,\cdot)=u_0$ and $|||u||| <+\infty$} \right\}.
$$
It is readily seen that $X$ is a complete metric space endowed
with the distance induced by the norm $|||\cdot|||$. For $u \in
X$, define the function
\begin{equation}\label{def: theta bonne}
\Theta u: t \in [0,T_\ast] \rightarrow K(t,\cdot) \ast u_0 -
\frac{1}{2} \int_{0}^{t}
\partial_x K (t-s,\cdot) \ast u^2(s,\cdot)ds \in L^2(\R).
\end{equation}
By Proposition \ref{prop: cont semi groupe}, $\Theta u \in
C([0,T_\ast];L^2(\R))$ and satisfies $\Theta u(0,\cdot)=u_0$. Define
$v:= u^2$. We have $\partial_x v=2 u \partial_x u$. Therefore that $v \in C([0,T_\ast];L^1(\R)) \cap C((0,T_\ast];W^{1,1}(\R))$
and that \eqref{hyp v reg} holds true. By Lemma \ref{lemme reg},
we deduce that $\Theta u \in X$. Let us take
$R>||u_0||_{L^2(\R)}+\mathcal{K}_1 ||u_0||_{L^2(\R)}$ and assume
that $|||u||| \leq R$. Since $||u^2||_{L^\infty((0,T_\ast);L^1(\R))} =||u||^2_{C([0,T_\ast];L^2(\R))},$
  estimate \eqref{esti: duhamel} of Remark
\ref{rem esti L 2} implies that
\begin{eqnarray}
||\Theta u||_{C([0,T_\ast];L^2(\R))} & \leq & e^{\omega_0 T_\ast}
||u_0||_{L^2(\R)}+ 2 \mathcal{K}_0 T_\ast^{\frac{1}{4}} \:
||u||^2_{C([0,T_\ast];L^2(\R))}, \nonumber \\
& \leq & e^{\omega_0 T_\ast} ||u_0||_{L^2(\R)}+ 2 \mathcal{K}_0
T_\ast^{\frac{1}{4}} \: R^2. \label{esti a reutiliser 1}
\end{eqnarray}
Estimate  \eqref{esti L2 gradient duhamel} of Lemma \ref{lemme
reg}, implies that
\begin{eqnarray*}
\sup_{t \in (0,T_\ast]} t^{\frac{1}{2}} ||\partial_x (\Theta
u(t,\cdot))||_{L^2(\R)} & \leq & \mathcal{K}_1 ||u_0||_{L^2(\R)} +
\frac{\mathcal{K}_0 I}{2} \; {T_\ast^\frac{1}{4}} \sup_{t \in
(0,T_\ast]} t^\frac{1}{2} ||\partial_{x}
(u^2) (t,\cdot)||_{L^1(\R)},\\
& \leq & \mathcal{K}_1 ||u_0||_{L^2(\R)} + \mathcal{K}_0 I \;
{T_\ast^\frac{1}{4}}  R^2,\\
\end{eqnarray*}
by Cauchy-Schwarz inequality.
Adding this inequality with \eqref{esti a reutiliser 1}, we get:
$$
|||\Theta u||| \leq e^{\omega_0 T_\ast}
||u_0||_{L^2(\R)}+\mathcal{K}_1 ||u_0||_{L^2(\R)} + \left(2+I
 \right)\mathcal{K}_0 T_\ast^{\frac{1}{4}}
R^2.
$$
For $T_\ast \in (0,T]$ sufficiently small such that
\begin{equation}\label{conditon duhamel 1}
e^{\omega_0 T_\ast} ||u_0||_{L^2(\R)}+\mathcal{K}_1
||u_0||_{L^2(\R)} + \left(2+I
 \right)\mathcal{K}_0 T_\ast^{\frac{1}{4}} R^2 \leq R,
\end{equation}
we deduce that $|||\Theta u||| \leq R$. To sum-up, we have
established that for any $T_\ast \in (0,T]$ such that
\eqref{conditon duhamel 1} holds true, $\Theta$ (defined by
\eqref{def: theta bonne}) maps $\overline{B}_R$ into
itself, where  $\overline{B}_R$ denotes the ball of $X$
 (endowed with the
$|||\cdot|||$ norm) centered at the origin and of radius $R$. Let
us now prove that $\Theta$ is a contraction. For $u,v \in
\overline{B}_R$, Estimate \eqref{esti: contraction bonne} of
Lemma \ref{rem: contraction} implies that
\begin{equation}\label{esti 3}
||\Theta u- \Theta v||_{C([0,T_\ast];L^2(\R))} \leq 4R
\mathcal{K}_0 T_\ast^{\frac{1}{4}} \:
||u-v||_{C([0,T_\ast];L^2(\R))},
\end{equation}
where we again used  $||u^2-v^2||_{C([0,T_\ast];L^1(\R))} \leq
(||u||_{C([0,T_\ast],L^2(\R))}+||v||_{C([0,T_\ast],L^2(\R))}) ||u-v||_{C([0,T_\ast];L^2(\R))}.$
Moreover, Estimate \eqref{contrac reg} of Lemma \ref{lemme reg}
implies that
\begin{eqnarray*}
\sup_{t \in (0,T_\ast]} t^{\frac{1}{2}}||\partial_x (\Theta u-
\Theta v)(t,\cdot)||_{L^2(\R)}  \leq  \mathcal{K}_0 I
T_\ast^{\frac{1}{4}} \sup_{t \in (0,T_\ast]} t^{\frac{1}{2}}||(u
\partial_xu -v \partial_x v)(t,\cdot)||_{L^1(\R)}.
\end{eqnarray*}
Since
\begin{eqnarray*}
t^{\frac{1}{2}} ||(u
\partial_xu -v \partial_x v)(t,\cdot)||_{L^1(\R)} & \leq &
t^{\frac{1}{2}} ||\partial_x v(t,\cdot)||_{L^2(\R)} ||
(u-v)(t,\cdot)||_{L^2(\R)}\\
& & \quad \quad +t^{\frac{1}{2}} ||u(t,\cdot)||_{L^2(\R)}
||\partial_x
(u-v)(t,\cdot)||_{L^2(\R)},\\
& \leq & |||v|||\; || (u-v)(t,\cdot)||_{L^2(\R)}\\
& & \quad \quad +|||u||| \;
t^{\frac{1}{2}} ||\partial_x (u-v)(t,\cdot)||_{L^2(\R)},\\
& \leq & R |||u-v|||,
\end{eqnarray*}
we get: $ \sup_{t \in (0,T_\ast]} t^{\frac{1}{2}}||\partial_x
(\Theta u- \Theta v)(t,\cdot)||_{L^2(\R)} \leq  R \mathcal{K}_0 I
T_\ast^{\frac{1}{4}} |||u-v|||. $ Adding this inequality with
\eqref{esti 3}, we find that
$$
|||\Theta u-\Theta v||| \leq (4+I)R \mathcal{K}_0
T_\ast^{\frac{1}{4}} |||u-v|||.
$$
Consequently, for any $T_\ast >0$ sufficiently small such that
\eqref{conditon duhamel 1} holds true and $(4+I)R \mathcal{K}_0
T_\ast^{\frac{1}{4}}<1$, $\Theta$ is a contraction from $
\overline{B}_R$ into itself. The Banach fixed point theorem then
implies that $\Theta$ admits a (unique) fixed point $u \in
C([0,T_\ast];L^2(\R)) \cap C((0,T_\ast];H^1(\R))$ satisfying
$\sup_{t \in (0,T_\ast]}
t^\frac{1}{2} ||\partial_{x} u(t,\cdot)||_{L^2(\R)} <\infty$ which is, of
course, a mild solution to \eqref{fowlereqn}.

To prove the $H^2$ regularity of $u$, we have to use again a
contracting fixed point theorem. But, this is now the gradient of
the solution which is searched as a fixed point. Let $t_0 \in
(0,T_\ast)$. For any $t \in (0,T_\ast-t_0]$, define
$\overline{u}(t,\cdot):=u(t_0+t,\cdot)$. Let $T_\ast' \in
(0,T_\ast-t_0]$. We still endow $C([0,T_\ast'];L^2(\R)) \cap
C((0,T_\ast'];H^1(\R))$ with the norm $|||\cdot|||$ defined in
\eqref{norme regulari} with $T_\ast$ replaced by $T_\ast'$. Define
the complete metric  space
$$
X':= \left\{v \in C([0,T_\ast'];L^2(\R)) \cap
C((0,T_\ast'];H^1(\R)) \mbox{ s.t. $v(0,\cdot)=v_0$ and $|||v|||
<+\infty$} \right\},
$$
where $v_0:=\partial_x \overline{u}(0,\cdot)$. For $v \in X'$,
define the function
\begin{equation}\label{def: gradient theta bonne}
\Theta' v: t \in [0,T_\ast'] \rightarrow K(t,\cdot) \ast v_0 -
\int_{0}^{t}
\partial_x K (t-s,\cdot) \ast (\overline{u}v)(s,\cdot)ds \in L^2(\R).
\end{equation}
Arguing as in the first step of the proof, we claim that
Proposition \ref{prop: cont semi groupe}, Remark \ref{rem esti L
2},   Lemmas \ref{rem: contraction} and \ref{lemme reg} imply that
$\Theta'$ maps $X'$ into itself with: for any $u,v \in X'$,
\begin{equation*}
|||\Theta' v||| \leq e^{\omega_0 T_\ast'}
||v_0||_{L^2(\R)}+\mathcal{K}_1 ||v_0||_{L^2(\R)} + C
T_\ast'^{\frac{1}{4}} |||v|||,
\end{equation*}
\begin{equation*}
|||\Theta' v-\Theta' w||| \leq C T_\ast'^{\frac{1}{4}} |||v-w|||,
\end{equation*}
for some nonnegative constant $C$ that only depends on
$\mathcal{K}_0$ and $||\overline{u}||_{C([t_0,T_\ast];H^1(\R))}$.
Let us take $R'$ such that
$$
R'
> e^{\omega_0 T_\ast'} ||v_0||_{L^2(\R)}+\mathcal{K}_1
||v_0||_{L^2(\R)}.
$$
If $T_\ast'>0$ satisfies
$$
e^{\omega_0 T_\ast'} ||v_0||_{L^2(\R)}+\mathcal{K}_1
||v_0||_{L^2(\R)} + C T_\ast'^{\frac{1}{4}} R' \leq R' \quad
\mbox{and} \quad C T_\ast'^{\frac{1}{4}} <1,
$$
then $\Theta'$ maps $\overline{B}_{R'}(X')$ into itself and is a
contraction. Let $v$ denote its unique fixed point. Observe now
that $\Theta' \partial_x \overline{u}=\partial_x \overline{u}$,
thanks to Remark \ref{rem for grad}. But, similar arguments than
these ones used to prove the uniqueness of a mild solution in the
preceding section allow to show that there exists at most one
function $w \in L^\infty((0,T_\ast');L^2(\R))$ that satisfies
$\Theta' w=w$. It follows that $\partial_x \overline{u}=v \in X'$
on $(0,T_\ast')$; hence, we deduce that $u \in
C(t_0,t_0+T_\ast'];H^2(\R))$. To sum-up, we have proved that for
all $t_0 \in (0,T_\ast]$, there exists $T_\ast' \in
(0,T_\ast-t_0]$ such that $u \in C((t_0,t_0+T_\ast'];H^2(\R))$. This
completes the proof of the continuity of $u$ on $(0,T_\ast]$ with
values in $H^2$. The proof of the $C^{1,2}$ regularity is
postponed to Lemma \ref{prop: reg bord} in the next section, where
it will be useful for the maximum principle failure.
\end{proof}

We can finally prove the global-in-time existence.
\begin{Prop}\label{prop: global}
Let $u_0 \in L^2(\R)$ and $T>0$. There exists a (unique) mild
solution $u \in C([0,T];L^2(\R)) \cap C((0,T];H^2(\R))$ to
\eqref{fowlereqn} such that
$$
\sup_{t \in (0,T]} t^\frac{1}{2} ||\partial_{x}
u(t,\cdot)||_{L^2(\R)} < +\infty \quad \mbox{and} \quad \sup_{t
\in (0,T]} t ||\partial_{xx}^2 u(t,\cdot)||_{L^2(\R)} < +\infty.
$$
Moreover, $u$ belongs to $C^{1,2}((0,T] \times
\R)$ and satisfies the PDE in \eqref{fowlereqn} in the classical
sense.
\end{Prop}

\begin{proof}
We have to derive first a $L^2$ estimate on the local regular
solution $u$ constructed in Proposition \ref{prop: exist local}.
Multiplying \eqref{fowlereqn} by $u$ and integrating w.r.t. the
space variable, we get:
\begin{equation}\label{tech esti L deux}
\frac{d}{dt} \frac{1}{2} \int_\R u^2dx+\int_\R
(\LL[u]-\partial_{xx}^2 u)u dx=0.
\end{equation}
Indeed, the following computations show that the nonlinear term
equals $0$:
$$
\int_\R \partial_x \left(\frac{u^2}{2}\right) u dx=-\int_\R
\frac{u^2}{2}
\partial_x u dx =-\frac{1}{2} \int_\R u (\partial_x u u)dx=-\frac{1}{2}\int_\R u \: \partial_x \left(\frac{u^2}{2}\right)
dx.
$$
But, Corollary \ref{proppseudodiff} implies that
\begin{eqnarray*}
\int_\R (\LL[u]-\partial_{xx}^2 u) u dx = \int_\R \F^{-1}(\psi_\LL
\F u) u dx = \int_\R \psi_\LL |\F u|^2 d\xi
 = \int_\R \mbox{Re}(\psi_\LL) |\F u|^2 d\xi,
\end{eqnarray*}
since $\int_\R (\LL[u]-\partial_{xx}^2 u) udx$ is real. It follows
that,
\begin{eqnarray*}
\int_{\R} (\LL[u]-\partial_{xx}^2 u) udx & \geq & \min
\mbox{Re}(\psi_\LL)  \int_\R |\F u|^2 d\xi, \\
& = & \min \mbox{Re}(\psi_\LL)  \int_\R u^2 dx,
\end{eqnarray*}
thanks to Plancherel's Equality. Equation \eqref{tech esti L deux}
then implies that
\begin{equation*}
\frac{d}{dt} \frac{1}{2} \int_\R u^2dx \leq \omega_0 \int_\R u^2dx
\end{equation*}
and by Gronwall's Lemma, we deduce that for all $t \in
[0,T_{\ast}]$
\begin{equation*}
||u(t,\cdot)||_{L^2(\R)} \leq e^{\omega_0 t} ||u_0||_{L^2(\R)}.
\end{equation*}
Define now
\begin{multline*}
t_0:=\sup \{t >0 \mbox{ s.t. there exists a (unique) mild
sol. to \eqref{fowlereqn}}\\ \mbox{on $(0,t)$ that satisfies the
regularity of Proposition \ref{prop: global}} \}
\end{multline*}
and let us assume that $t_0 < T$ (recall that Proposition
\ref{prop: exist local} ensures that $t_0>0$). By  Proposition
\ref{prop: exist local}, there exists $T_\ast >0$ such that for
any initial data $v_0$ that satisfy $||v_0||_{L^2(\R)} \leq
e^{\omega_0 t_0} ||u_0||_{L^2(\R)}$, \eqref{fowlereqn} admits a
regular mild solution on $(0,T_\ast)$ with initial datum $v_0$.
Hence, if we define $v_0=u(t_0-T_\ast/2)$, then \eqref{fowlereqn}
admits a mild solution $v$ that satisfies the regularity of
Proposition \ref{prop: global}. Using the uniqueness and the
semi-group property, it is now easy to show that $u(t_0
-T_\ast/2\, +t,\cdot)=v(t,\cdot)$ for all $t \in [0,T_\ast/2]$ and
that the function $\widetilde{u}$ defined by $\widetilde{u}=u$ on
$[0,t_0]$ and $\widetilde{u}(t_0-T_\ast/2\, +t,\cdot)=v(t,\cdot)$ for
$t \in [T_\ast/2,T_\ast]$ is still a mild solution to
\eqref{fowlereqn} that satisfies the regularity of proposition
\ref{prop: global}. Since the solution $\widetilde{u}$ lives on
$[0,t_0+T_\ast/2]$, this gives us a contradiction. We conclude that
$t_0 \geq T$ and this completes the proof of the global existence
of a regular solution.
\end{proof}
\begin{remark}\label{rem reg}
To sum-up, we have proved Theorem \ref{theo principal} with the
$C^{1,2}$ regularity of $u$. To obtain further regularity, we
claim that we can use the same method by arguing by induction.
\end{remark}
Now we  prove the $L^2$-stability stated in  Proposition \ref{L2stab}.
\begin{proof}[Proof of Proposition \ref{L2stab} ]
Let $(u,v)$ be solutions to \eqref{fowlereqn} with respective $L^2$ initial data $(u_0,v_0)$.
let $T>0$ and $t\in [0,T]$.
 Substracting
$$ u(t,\cdot)= K(t,\cdot) \ast u_0 -
\frac{1}{2} \int_{0}^{t}
\partial_x K (t-s,\cdot) \ast u^2(s,\cdot)ds $$
and
$$ v(t,.)= K(t,\cdot) \ast v_0 -
\frac{1}{2} \int_{0}^{t}
\partial_x K (t-s,\cdot) \ast v^2(s,\cdot)ds
$$
we get
\begin{equation}
  \label{soludiff}
  u(t,\cdot)- v(t,\cdot) = K(t,\cdot) \ast (u_0 - v_0) -
\frac{1}{2} \int_{0}^{t}
\partial_x K (t-s,\cdot) \ast (u^2(s,\cdot)- v^2(s,\cdot))ds.
\end{equation}
Hence, by \eqref{esti omega 0} of Remark \ref{omega 0} and Young inequality
$$|| u(t,\cdot) - v(t,\cdot)||_{L^2 (\R)} \leq e^{\omega_0 T} ||u_0 - v_0||_{L^2 (\R)} +
\frac{1}{2}\int_0^t ||\partial_x K (t-s,\cdot) ||_{L^2 (\R)}|| u^2(s,\cdot) - v^2(s,\cdot)||_{L^1 (\R)}\,ds.
$$
Taking $M = \max \left( ||u||_{C([0,T];L^2(\R))},
||v||_{C([0,T];L^2(\R))}\right), $ we can bound
\begin{eqnarray*}
|| u(t,\cdot) - v(t,\cdot)||_{L^2 (\R)}& \leq & e^{\omega_0 T}
||u_0 - v_0||_{L^2 (\R)} + M \int_0^t ||\partial_x K (t-s,\cdot)
||_{L^2
(\R)}|| u(s,\cdot) - v(s,\cdot)||_{L^2 (\R)}\,ds \\
& \leq &e^{\omega_0 T} ||u_0 - v_0||_{L^2 (\R)} + M \mathcal{K}_0
\int_0^t (t-s)^{-\frac{3}{4}}|| u(s,\cdot) - v(s,\cdot)||_{L^2
(\R)}\,ds,
\end{eqnarray*}
thanks to \eqref{esti gradient kernel}.
With lemma \ref{supergronwall}, the proof is finished.
\end{proof}

\section{Failure of the maximum principle}
\label{sect princ max}

We now investigate the proof of Theorem \ref{theo max}. We first
need a regularity result which ensures that if the initial data is
regular then so is the solution up to the time $t=0$.
\begin{lemme}\label{prop: reg bord}
Let $u_0 \in H^2(\R)$ and $T>0$. Assume that $u$ is a mild
solution to \eqref{fowlereqn} that satisfies the regularity of
Proposition \ref{prop: exist local}. Then, $u$ is in fact
$C([0,T];H^2(\R)) \cap C^{1,2}((0,T] \times \R)$ and satisfies the
PDE in \eqref{fowlereqn} in the classical sense. Moreover, if $u_0
\in C^2(\R)$, then $u \in C^{1,2}([0,T] \times \R)$ and satisfies
the PDE up to the time $t=0$.
\end{lemme}

\begin{proof}
First, we leave it to the reader to verify that the continuity
with values in $H^2$ up to the time $t=0$ can be proved again by
the use of a contracting fixed point theorem. Note that the
regularity of $u_0$ allows to work in a space of continuous
functions with values in $H^2$ up to the time $t=0$; more
precisely, we argue as in the proof of Proposition \ref{prop:
exist local}, but we can directly use the $C([0,T_\ast];H^2)$ norm
instead of the $|||\cdot|||$ norm defined in \eqref{norme
regulari}. Let us now prove that $u$ is a classical solution to
\eqref{fowlereqn}. Taking the Fourier transform w.r.t. the space
variable in \eqref{formule duhamel}, we get: for all $t \in
[0,T]$,
\begin{equation}\label{eq fourier}
\F (u(t,\cdot)) = e^{-t \psi_\LL} \F u_0-\int_0^t i \: \pi \cdot
e^{-(t-s) \psi_\LL} \F (u^2(s,\cdot)) ds.
\end{equation}
Since $u^2 \in C([0,T];L^1(\R))$, we know that $\F (u^2) \in
C([0,T];C_b(\R,\mathbb{C}))$. For any $\xi \in \R$, the function
$t \in [0,T] \rightarrow \F (u^2(t,\cdot))(\xi) \in \mathbb{C}$ is
thus continuous. Define
$$
w(t,\xi):=-\int_0^t i \: \pi \xi e^{-(t-s) \psi_\LL(\xi)} \F
(u^2(s,\cdot))(\xi)ds.
$$
Classical results on ODE then imply that $w$ is derivable w.r.t.
the time variable with
\begin{equation}\label{tech regularity}
\partial_t w(t,\xi)+\psi_\LL(\xi) w(t,\xi)= - i \: \pi \xi \F (u^2(t,\cdot))(\xi)=-\F\left(\partial_x \left(\frac{u^2}{2} \right)(t,\cdot) \right)(\xi).
\end{equation}
Let us prove that all these terms are continuous with values in
$L^2$. First, $u \in C([0,T];H^1(\R))$ therefore
$\partial_x (u^2) \in C([0,T];L^2(\R))$  and we deduce that
$\F(\partial_x(\frac{u^2}{2}))$ is continuous with values in
$L^2$. Moreover, Equation \eqref{eq fourier} implies that
$$
\psi_\LL \; w(t,\cdot)= \psi_\LL \left(\F (u(t,\cdot))-e^{-t
\psi_\LL} \F u_0 \right).
$$
Since $u \in C([0,T];L^2(\R)) $ and $\psi_\LL$ behaves at infinity
as $|\cdot|^2$, $\psi_\LL w$ is continuous with values in $L^2$.
All the terms in \eqref{eq fourier} then are continuous with
values in $L^2$ and this implies $w \in
C^1([0,T];L^2(\R,\mathbb{C}))$ with
$$
\frac{d}{dt} (w(t,\cdot))+\psi_\LL \; w(t,\cdot)=
-\F\left(\partial_x \left(\frac{u^2}{2} \right)(t,\cdot) \right).
$$
Moreover, it is easy to see that $t \in [0,T] \rightarrow e^{-t
\psi_\LL} \F u_0 \in L^2(\R,\mathbb{C})$ is $C^1$ with
$$
\frac{d}{dt} \left(e^{-t \psi_\LL} \F u_0 \right)+\psi_\LL e^{-t
\psi_\LL} \F u_0=0.
$$
From Equation \eqref{eq fourier}, we   infer that $\F u$ is $C^1$
on $[0,T]$ with values in $L^2$ with
\begin{multline*}
\frac{d}{dt} (\F( u(t,\cdot)))=-\psi_\LL w(t,\cdot)-\psi_\LL e^{-t \psi_\LL} \F u_0-\F\left(\partial_x \left(\frac{u^2}{2} \right)(t,\cdot) \right)\\
=-\psi_\LL \F(u(t,\cdot))-\F\left((u \partial_x u)(t,\cdot)
\right).
\end{multline*}
Since $\F$ is an isometry of $L^2$, we deduce that $u \in
C^1([0,T];L^2(\R))$ and that
\begin{eqnarray*}
\frac{d}{dt} (u(t,.)) & = & -\partial_x \left(\frac{u^2}{2} \right)(t,\cdot)- \F^{-1}\left(\psi_\LL \F(u(t,\cdot)) \right) ,\\
& = &  -\partial_x \left(\frac{u^2}{2}
\right)(t,\cdot)-\LL[u(t,\cdot)]+\partial_{xx}^2 u(t,\cdot),
\end{eqnarray*}
where we used Corollary \ref{proppseudodiff} to compute the
pseudo-differential term. In particular, $u$ satisfies the PDE of
\eqref{fowlereqn} in the distribution sense. What is left to prove
is the $C^2$ regularity in space of $u$. Differentiating
\eqref{formule duhamel} two times w.r.t. the space variable, we
get: for any $t \in [0,T]$,
\begin{equation}\label{formule: duhamel hessian}
\partial^2_{xx}  u (t,\cdot)= K(t,\cdot) \ast u_0'' -
\int_{0}^{t}
\partial_x K (t-s,\cdot) \ast v(s,\cdot) ds.
\end{equation}
where $v=(\partial_x u)^2+u
\partial_{xx}^2u$. By the Sobolev imbedding $H^2(\R) \hookrightarrow
C^1_b(\R)$, we know that $v \in C([0,T];L^1(\R) \cap L^2(\R))$. By
Lemma \ref{dernier}, we know that for all $x,y \in \R$,
$$
|\partial_x K (t-s,\cdot) \ast v(s,\cdot)(x)-\partial_x K
(t-s,\cdot) \ast v(s,\cdot)(y)| \leq ||\partial_x
K(t-s)||_{L^2(\R)} ||\mathcal{T}_{(x-y)}(v(s,\cdot))- v(s,\cdot) ||_{L^2(\R)}.
$$
By \eqref{esti gradient kernel}, we deduce that for all $t \in
[0,T]$ and all $x,y \in \R$,
\begin{multline*}
\left|\int_{0}^{t}
\partial_x K (t-s,\cdot) \ast v(s,\cdot)(x) ds-\int_{0}^{t}
\partial_x K (t-s,\cdot) \ast v(s,\cdot)(y) ds\right| \\
\leq \int_0^t \mathcal{K}_0 (t-s)^{-\frac{3}{4}}
||\mathcal{T}_{(x-y)}(v(s,\cdot))-v(s,\cdot)||_{L^2(\R)} ds \leq
4T^{\frac{1}{4}} \sup_{s \in [0,T]}
||\mathcal{T}_{(x-y)}(v(s,\cdot))-v(s,\cdot)||_{L^2(\R)}.
\end{multline*}
By Lemma \ref{lem tech 1}, we deduce that the second term of
\eqref{formule: duhamel hessian} is continuous w.r.t. the space
variable independently of the time variable (equicontinuity w.r.t. the time variable). Moreover, we already
know that this term is continuous on $[0,T]$ with values in $L^2$ (by Proposition \ref{prop: cont semi groupe})
and Lemma \ref{lemme tech 2} implies that it is continuous w.r.t.
the couple $(t,x)$ on $[0,T] \times \R$. We now leave it to the
reader to verify that $(t,x) \rightarrow K(t,\cdot) \ast u_0''(x)$
is continuous on $(0,T] \times \R$ when $u_0 \in H^2(\R)$ and
continuous on $[0,T] \times \R$ when moreover $u_0 \in C^2(\R)$.
The proof of Lemma \ref{prop: reg bord} is complete.
\end{proof}

The proof of Theorem \ref{theo max} is now an immediate
consequence of the integral formula \eqref{intformula}.
\begin{proof}[Proof of Theorem \ref{theo max}]
Lemma \ref{prop: reg bord} and Proposition
\ref{propintegralformula} imply that the solution $u$ to
\eqref{fowlereqn} is $C^{1,2}$ up to the initial time $t=0$ and
that
$$
u_t(0,x_\ast)+u_0(x_\ast) u_0'(x_\ast)+C_{\LL} \int_{-\infty}^{0}
\frac{u_0(x_\ast+z)-u_0(x_\ast)-u_0'(x_\ast) z}{|z|^{7/3}} \:
dz-u_0''(x_\ast)=0.
$$
It follows that
$$
u_t(0,x_\ast)=-C_{\LL} \int_{-\infty}^{0}
\frac{u_0(x_\ast+z)}{|z|^{7/3}} \: dz<0.
$$
There then exists $t_\ast >0$ such that $u(t_\ast,x_\ast)<0$. The
proof of Theorem \ref{theo max} is now complete.
\end{proof}

\section{Numerical simulations}
\label{sect numerique}

The aim of this part is to show some numerical simulations for
\eqref{fowlereqn}. An explicit discretization
gives results in line with the theoretical study (see Remark \ref{rem max}).

We write \eqref{fowlereqn} with a viscous coefficient
$\varepsilon>0$ as follows:
\begin{equation}
\partial_t u + \partial_x \left(\frac{u^2}{2} + \LLL [u]\right)- \varepsilon \partial^2_{xx} u = 0,
\label{fowlereqn_num}
\end{equation}
where for any $\varphi \in \mathcal{S}(\R)$ and $x \in \R$,
\begin{equation*}
\LLL[\varphi](x) := \int_{0}^{+\infty} |\zeta|^{-\frac{1}{3}}
\varphi' (x-\zeta) d\zeta. \label{nonlocalterm_num}
\end{equation*}
The viscous coefficient is taken sufficiently small, in order to
magnify the erosive effect of the non-local term. The new
definition of the non-local term ($\LL[u]=\partial_x \LLL[u]$)
follows \cite{Fow01}, which interpretes $\LLL[u]$ as a flow.
Notice that in \cite{Fow01,Fow02}, the bottom is, in fact,
$s(t,x)=u(t,x+q'(1) t)$, where $q$ is the bedload transport of
sediments; for the sake of simplicity, we continue to work with
$u$.

To shed light on the effect of the nonlocal term, we compare the
evolution of the solution of \eqref{fowlereqn_num} with the solution
of the viscous Burgers equation:
\begin{equation}
\partial_t u + \partial_x \left(\frac{u^2}{2}\right)- \varepsilon \partial^2_{xx} u = 0.
\label{burgers}
\end{equation}

\subsection{Maximum principle for the viscous Burgers equation}

\begin{figure}[!ht!]
\begin{center}
\includegraphics[scale=0.325]{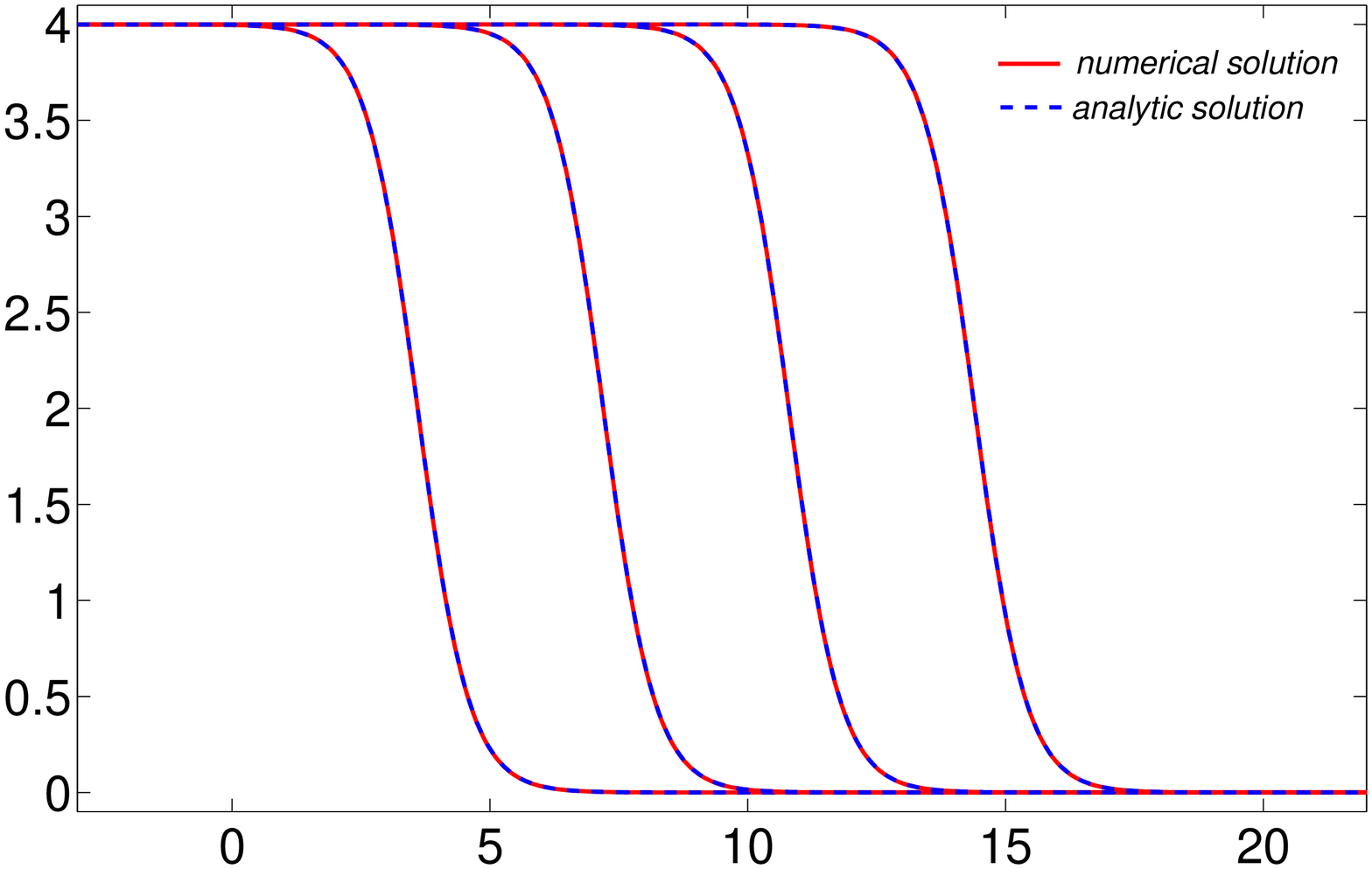}
\caption{Numerical and analytic travelling waves of the viscous
Burgers equation.} \label{sol_anal1}
\end{center}
\end{figure}

\begin{figure}[!ht!]
\begin{center}
\includegraphics[scale=0.5]{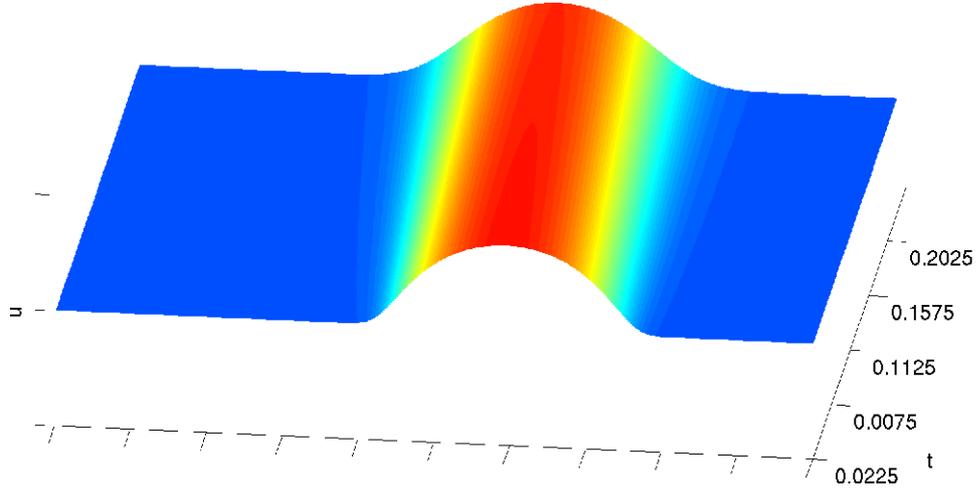}
\caption{Evolution of the solution of \eqref{burgers} with $u_0$
defined in \eqref{initial_dune} ($L=30$, $M=4001$ and
$\varepsilon=0.1$). \label{num_burg1}}
\end{center}
\end{figure}

It is well-known that \eqref{burgers} satisfies the maximum
principle: for any initial data $u_0 \in L^\infty(\R)$,
$\mbox{ess-inf} \: u_0  \leq u \leq \mbox{ess-sup} \: u_0$. As a
consequence, \eqref{burgers} cannot
take into account erosion phenomena. To simulate the evolution of
$u$, we define a regular discretization of $[0,L]$ with a spatial
step $\Delta x$ such that $L=M\Delta x$, and a discretization of
$[0,T]$ with a time step $\Delta t$ such that $T=N \Delta t$. We
let $x_i$, $t_n$ and $u_i^n$ respectively denote the point
$i\Delta x$, the time $n\Delta t$ and the computed solution at the
point $(n \Delta t,i \Delta x)$. We use the following explicit
centered scheme:
\begin{equation}
u_i^{n+1}  =  u_i^n+\Delta t
\left[-\frac{1}{2}\frac{{(u_{i+1}^{n})}^{2}-{(u_{i-1}^{n})}^{2}}{2\Delta
x} + \varepsilon \frac{u_{i+1}^{n}-2u_{i}^{n}+u_{i-1}^{n}}{{\Delta
x}^2} \right]. \label{shemaBURG_explicit}
\end{equation}
It is well-known that this scheme is stable under the CFL-Peclet
condition:
\begin{equation}
\Delta t = \min \left( \frac{\Delta x}{|u|} ,\frac{{\Delta x}^2}{2
\varepsilon} \right). \label{stab_cond}
\end{equation}
To convince the reader, let us simulate the evolution of the
well-known following travelling waves of \eqref{burgers} for
$\varepsilon=1$:
\begin{equation*}
u(t,x):=\frac{1}{2} \left[ 1-\tanh \left( \frac{1}{4}\left(x-\frac{1}{2} t\right) \right) \right].
\end{equation*}
We expose in Figure \ref{sol_anal1} both analytic and numerical
solutions. We observe an error of the order of $10^{-4}$  between these solutions.
Let us now take, as an initial dune, the following small regular
perturbation on the bottom:
\begin{equation}
u_0(x)= \left\{
\begin{array}{ll}
e^{\frac{-1}{1-(x-\frac{L}{2})^2}} & \mbox{if } \frac{L}{2}-1 < x
< \frac{L}{2}+1, \\
0 & \mbox{otherwise}.
\end{array}
\right.\label{initial_dune}
\end{equation}
We describe its evolution in Figure \ref{num_burg1}. The dune
propagates, but as mentioned above the erosion phenomena are not
taken into account since $u$ remains positive (because of the
maximum principle).
\begin{remark}
    Equation \eqref{fowlereqn} also admits travelling wave solutions,  see \cite{borysaz}.
\end{remark}

\subsection{Erosive effect of the nonlocal term}

\begin{figure}[!ht!]
\begin{center}
\includegraphics[height=4cm,width=8cm]{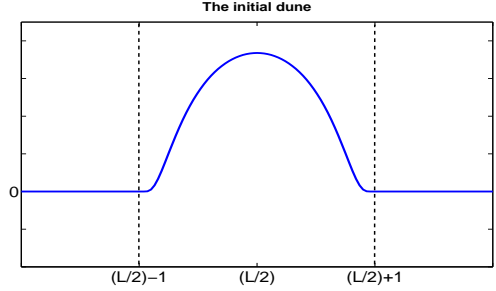}
\caption{The initial dune defined in \eqref{initial_dune}.}
\label{num_burg2}
\end{center}
\end{figure}

Let us return to the study of \eqref{fowlereqn_num}. We add the
discretization of the non-local operator $\LLL$ to the explicit
centered scheme (\ref{shemaBURG_explicit}). It is natural to
consider the following discretization:
\begin{equation*}
\LLL[u_i^n] \approx \sum_{j=0}^{+\infty} |j
\Delta x|^{-\frac{1}{3}}
\frac{u_{i-j+1}^{n}-u_{i-j-1}^{n}}{2{\Delta x}}, \label{L_discret}
\end{equation*}
Instead we will approximate
\begin{equation*}
\LLL[u_i^n] \approx \sum_{j=0}^{i} |j
\Delta x|^{-\frac{1}{3}}
\frac{u_{i-j+1}^{n}-u_{i-j-1}^{n}}{2{\Delta x}}. \label{L_discret}
\end{equation*}
This  is based on the assumption that for $x \in [0,L]$,
\begin{equation}\label{discret_nonlocal}
\LLL[u(t,\cdot)](x) \approx \int_{0}^{x} |\zeta|^{-\frac{1}{3}}
\partial_x u (t,x-\zeta) d\zeta.
\end{equation}
This fact is not true for general $u$, but if we
assume that the  initial profile $u_0$
satisfies
$u_0(x) = 0, \forall x  \leq 0$ and semi-discretize in time
Equation \eqref{fowlereqn_num}, we get :
$$ u(t+\Delta t,x) = u(t,x) + \Delta t \left( -\partial_x (\frac{u^2}{2}) -\partial_x \LLL [u(t,.)] + \varepsilon \partial_{xx}^2 u
\right).
$$
We observe that $u(t+\Delta t,x) = 0, \;\forall x  \leq 0$ and by
induction $u(t_n, x)= 0   \;\forall x  \leq 0,\; \forall n$. Now
$$\LLL[u(t_n,.)] = \int_{0}^{x} |\zeta|^{-\frac{1}{3}}
\partial_x u (t_n,x-\zeta) d\zeta.$$
Actually, we take  $u_0 \in C_c^\infty(\R)$ and $\mbox{supp}(u_0) \subset
\subset (0,L)$ (see Figure \ref{num_burg2}).
Moreover, Lemma \ref{prop: reg bord} suggests that all the
derivatives of $u$ are continuous with values in $L^2$ w.r.t. the
time variable up to the time $t=0$. It then is natural to expect
that (at least for small times) Equation \eqref{discret_nonlocal}
is a good approximation.

\begin{figure}[!ht!]
\begin{center}
\includegraphics[scale=0.5]{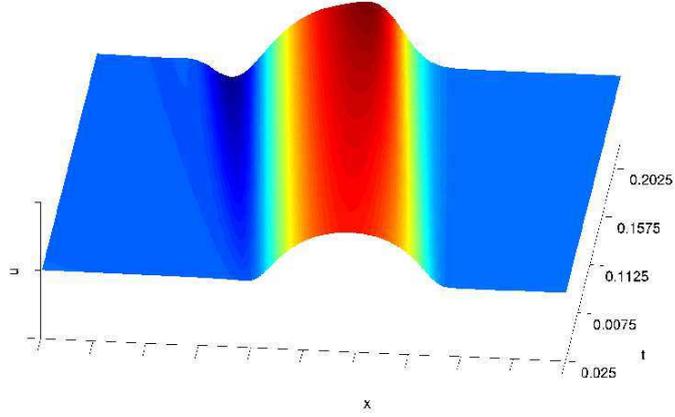}
\includegraphics[scale=0.8]{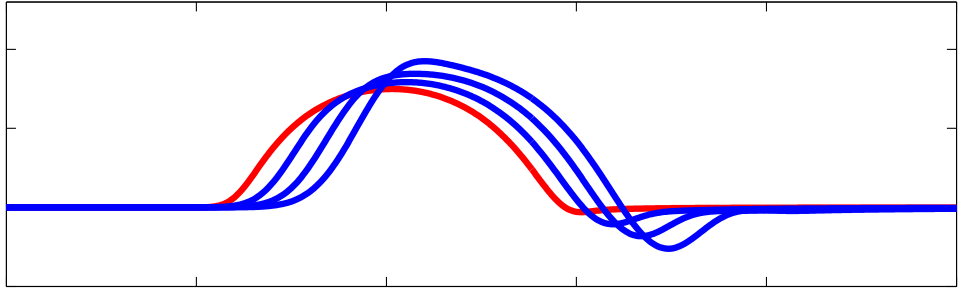}
\caption{Evolution of an initial dune, by using the non-local
model (\ref{fowlereqn_num}) ($L=30$, $M=4001$ and
$\varepsilon=0.1$).} \label{num_burg3}
\end{center}
\end{figure}

We then use the following explicit scheme for (\ref{fowlereqn_num}):
\begin{equation*}
u_i^{n+1}  =  u_i^n+\Delta t
\left[-\frac{1}{2}\frac{{(u_{i+1}^{n})}^{2}-{(u_{i-1}^{n})}^{2}}{2\Delta
x} -\frac{\LLL[u_{i+1}^n]-\LLL[u_{i-1}^n]}{2\Delta x} +
\varepsilon \frac{u_{i+1}^{n}-2u_{i}^{n}+u_{i-1}^{n}}{{\Delta
x}^2} \right]. \label{shemaFOW_explicit}
\end{equation*}
As far as the stability condition, one can numerically see that
(\ref{stab_cond}) is still ensuring stability for small $\Delta x$. The
evolution of the initial dune (\ref{initial_dune}) is given in
Picture \ref{num_burg3}. As the solutions of the viscous Burgers
equation, the dune is propagated downstream  but we now observe an erosive
process behind the dune: the bottom is eroded downstream from the
dune, as shown in Remark \ref{rem max}. \\
Let us make a final remark.
We are aware of  that the fact that these numerical simulations are a first crude attempt.
To tackle rigorously the non local term would need further study, which will be reported elsewhere.

\newpage

\appendix

\section{Some technical lemmas.}
 We first recall  a generalization of Gronwall's lemma proved e.g. in \cite{fontbona}.
\begin{lemme}\label{supergronwall}
 Let $g : [0,T] \rightarrow \R_+$ be a bounded measurable function and suppose that there are positive constants $C, A$ and $\theta >0$ such that, for all $t\leq T$ ,
$$ g(t) \leq A + C \int_0^t (t-s)^{\theta -1} g(s)\,ds.
$$
Then,
$$ \sup_{0\leq t\leq T} g(t) \leq C_T A,$$
where constant $C_T$ does not depend on $A$.
\end{lemme}

\begin{lemme}\label{dernier}
Let $f,g \in L^2(\R)$. Then, $f \ast g \in C(\R)$ and for all $x,y
\in \R$,
$$
|f \ast g(x)-f \ast g(y)| \leq ||\mathcal{T}_{(x-y)}
f-f||_{L^2(\R)} ||g||_{L^2(\R)}.
$$
\end{lemme}
\begin{proof}
The result is immediate if $f$ and $g$ are smooth; indeed,
\begin{eqnarray*}
|f \ast g(x)-f \ast g(y)|& =&\left|\int_\R f(x-z) g(z) dz-\int_\R
f(y-z) g(z) dz \right|,\\
& \leq & \int_\R |f(x-z)-f(y-z) g(z) |dz,\\
& \leq & ||\mathcal{T}_{(x-y)} f-f||_{L^2(\R)} ||g||_{L^2(\R)}.
\end{eqnarray*}
The result for general $f$ and $g$ only $L^2$, is then obtained by
density.
\end{proof}

\begin{lemme}\label{lem tech 1}
Let $u \in C([0,T];L^2(\R))$. Then, $ \sup_{t \in [0,T]}
||\mathcal{T}_{h} (u(t,\cdot)) -u(t,\cdot)||_{L^2(\R)} \rightarrow
0, $ as $h \rightarrow 0$.
\end{lemme}

\begin{proof}
The function $u$ is uniformly continuous with values in $L^2$ as a
continuous function on a compact set $[0,T]$. For any
$\varepsilon>0$, there then exist finite a sequence $0=t_0 <t_1<
\ldots <t_N=T$ such that for any $t \in [0,T]$, there exists $j
\in \{0,\ldots,N-1\}$ with
$$
||u(t,\cdot)-u(t_j,\cdot)||_{L^2(\R)} \leq \varepsilon.
$$
Moreover,
\begin{multline*}
||\mathcal{T}_{h} (u(t,\cdot)) -u(t,\cdot)||_{L^2(\R)} \leq
||\mathcal{T}_{h} (u(t,\cdot))
-\mathcal{T}_{h}(u(t_j,\cdot))||_{L^2(\R)}\\
+||\mathcal{T}_{h}(u(t_j,\cdot))-
u(t_j,\cdot)||_{L^2(\R)}+||u(t_j,\cdot) -u(t,\cdot)||_{L^2(\R)}.
\end{multline*}
Since $||\mathcal{T}_{h} (u(t,\cdot))
-\mathcal{T}_{h}(u(t_j,\cdot))||_{L^2(\R)}=||u(t,\cdot)
-u(t_j,\cdot)||_{L^2(\R)}$, we get:
\begin{eqnarray*}
||\mathcal{T}_{h} (u(t,\cdot)) -u(t,\cdot)||_{L^2(\R)} & \leq &
||\mathcal{T}_{h}(u(t_j,\cdot))-
u(t_j,\cdot)||_{L^2(\R)}+2||u(t_j,\cdot) -u(t,\cdot)||_{L^2(\R)},\\
& \leq & ||\mathcal{T}_{h}(u(t_j,\cdot))-
u(t_j,\cdot)||_{L^2(\R)}+2 \varepsilon.
\end{eqnarray*}
By the continuity of the translation in $L^2(\R)$,
$||\mathcal{T}_{h}(u(t_j,\cdot))- u(t_j,\cdot)||_{L^2(\R)}
\rightarrow 0$, as $h \rightarrow 0$. Then,
$$
\limsup_{h \rightarrow 0} ||\mathcal{T}_{h} (u(t,\cdot))
-u(t,\cdot)||_{L^2(\R)} \leq 2 \varepsilon.
$$
Taking the infimum w.r.t. $\varepsilon>0$ implies the result.
\end{proof}

\begin{lemme}\label{lemme tech 2}
Let $u \in C([0,T];L^2(\R))$ such that $u$ is continuous w.r.t. the variable
$x$ uniformly in $t$. Then, $u \in C([0,T] \times \R)$.
\end{lemme}

\begin{proof}
Let $(t_0,x_0) \in [0,T] \times \R$. Let $ \varepsilon >0$. By the
regularity of $u$ w.r.t. the space variable, we know that there
exists $\eta>0$ such that for any $t \in [0,T]$ and all $x,y \in
[x_0-\eta,x_0+\eta]$,
\begin{eqnarray*}
|u(t_0,x_0)-u(t,x)| & \leq & |u(t_0,x_0)-u(t_0,y)|+|u(t_0,y)-u(t,y)|+|u(t,y)-u(t,x)|,\\
& \leq & \varepsilon+ |u(t_0,y)-u(t,y)|+\varepsilon.
\end{eqnarray*}
If we integrate w.r.t. $y \in [x_0-\eta,x_0+\eta]$, then we get:
$$
2 \eta |u(t_0,x_0)-u(t,x)| \leq 4 \varepsilon \eta
+\int_{x_0-\eta}^{x_0+\eta} |u(t_0,y)-u(t,y)|dy \leq 4 \varepsilon
\eta (2\eta)^{\frac{1}{2}} ||u(t_0,\cdot)-u(t,\cdot)||_{L^2(\R)}.
$$
By the continuity of $u$ with values in $L^2$,
$$
\limsup_{(t,x) \rightarrow (t_0,x_0)} |u(t_0,x_0)-u(t,x)| \leq 2
\varepsilon.
$$
Taking the infimum w.r.t. $\varepsilon>0$ completes the proof.
\end{proof}

\textbf{Acknowledgements.} We thank B. Mohammadi for advice on the
numerical scheme and B. Alvarez-Samaniego for helpful comments. P. Azerad and D. Is\`ebe are supported by the ANR
(project COPTER NT05-2-42253).

\end{document}